\def\year{2022}\relax
\documentclass[letterpaper]{article}
\usepackage{aaai22}
\usepackage{times}  % DO NOT CHANGE THIS
\usepackage{helvet}  % DO NOT CHANGE THIS
\usepackage{courier}  % DO NOT CHANGE THIS
\usepackage[hyphens]{url}  % DO NOT CHANGE THIS
\usepackage{graphicx} % DO NOT CHANGE THIS
\usepackage{natbib}  % DO NOT CHANGE THIS AND DO NOT ADD ANY OPTIONS TO IT
\usepackage{caption} % DO NOT CHANGE THIS AND DO NOT ADD ANY OPTIONS TO IT
\DeclareCaptionStyle{ruled}{labelfont=normalfont,labelsep=colon,strut=off} % DO NOT CHANGE THIS
\usepackage{amsmath}
\usepackage{amsfonts}
\usepackage{multirow}
\usepackage{mathtools}
\usepackage{amsthm}
\usepackage{bm}
\usepackage{makecell}
\usepackage{xcolor}

\newcommand{\norm}[1]{\left\lVert#1\right\rVert}
\newcommand{\bmm}{m^*}
\newcommand{\bmmh}{m_h^*}
\newcommand{\romanone}{\uppercase\expandafter{\romannumeral1}}
\newcommand{\romantwo}{\uppercase\expandafter{\romannumeral2}}
\newtheorem{prop}{Proposition}

\title{Solving PDE-Constrained Control Problems Using Operator Learning}
\author{
    Rakhoon Hwang\equalcontrib\thanks{Work performed while at Pohang University of Science and
Technology},\textsuperscript{\rm 1} 
    Jae Yong Lee\equalcontrib\footnotemark[2],\textsuperscript{\rm 2}
    Jin Young Shin\equalcontrib,\textsuperscript{\rm 3} 
    Hyung Ju Hwang\thanks{Corresponding author.}\textsuperscript{\rm 3}
}
\affiliations{
    \textsuperscript{\rm 1}POSTECH Institute of Artifcial Intelligence\\
    \textsuperscript{\rm 2}Center for Artifcial Intelligence and Natural Sciences, Korea Institute for Advanced Study\\
    \textsuperscript{\rm 3}Department of Mathematics, Pohang University of Science and Technology\\
    hrh0125@postech.ac.kr, jaeyong@kias.re.kr, \{sjy6006, hjhwang\}@postech.ac.kr
}
\usepackage{lipsum}

\begin{document}

\maketitle

\begin{abstract}
The modeling and control of complex physical systems are essential in real-world problems. We propose a novel framework that is generally applicable to solving PDE-constrained optimal control problems by introducing surrogate models for PDE solution operators with special regularizers. The procedure of the proposed framework is divided into two phases: solution operator learning for PDE constraints (Phase 1) and searching for optimal control (Phase 2). Once the surrogate model is trained in Phase 1, the optimal control can be inferred in Phase 2 without intensive computations. Our framework can be applied to both data-driven and data-free cases. We demonstrate the successful application of our method to various optimal control problems for different control variables with diverse PDE constraints from the Poisson equation to Burgers' equation.
\end{abstract}

\section*{Introduction} \label{intro}
\noindent The modeling of physical systems to support decision making and solve the optimal control problem is a key problem in many industrial, economical, and medical applications. Such systems can be described mathematically through partial differential equations (PDEs). In this regard, solving PDE-constrained optimal control problems provides a control law for a complex system governed by PDEs. A PDE-constrained optimal control problem has been successfully used in many applications: shape optimization \cite{haslinger2003introduction,sokolowski1992introduction}, mathematical finance \cite{bouchouev1999uniqueness,egger2005tikhonov}, and flow control \cite{gunzburger2002perspectives}. 

As computing power increases and optimization technologies significantly improve, many researchers have studied algorithms and computational methods that are accurate, efficient, and applicable for complex physical systems. In control theory, adjoint methods are some of the most common approaches to handle this problem. Adjoint methods provide an efficient way to compute gradients that appear in optimization problems \citep{lions1971optimal, pironneau1974optimum, troltzsch2010optimal}, and its variants have been applied to many different areas of science. However, adjoint-based iterative schemes, such as shooting methods, suffer from computational costs because of their trial-and-error nature. In addition, iterative schemes are often sensitive to the initial guesses to the solutions.

Deep learning methods have recently derived major techniques for scientific computations, including PDEs or optimization problems \cite{raissi2019physics}. In particular, solving a family of parametrized PDEs requires networks to approximate a function-to-function mapping or operator. In \cite{li2020fourier}, \cite{li2020neural}, \cite{lu2019deeponet}, \cite{ zhu2018bayesian}, and \cite{zhu2019physics}, the authors utilized neural networks to learn a mapping from the parameters (e.g. initial or boundary) of a PDE to the corresponding solution. Such models are used as a surrogate model to solve various problems such as uncertainty quantification \cite{zhu2018bayesian,zhu2019physics} and an inverse problem \cite{li2020fourier}.

Recent studies have been conducted to solve PDE-constrained optimal control problems through a deep learning approach. In \cite{holl2020learning}, the authors proposed a predictor-corrector scheme for long-term fluid dynamics control, combining neural networks with a differentiable solver. In \cite{rabault2019artificial}, the authors experimentally showed that active flow control for vortex shedding and a drag reduction can be achieved through model-free reinforcement learning. Model-free methods are often known to require numerous interactions with the environment to search for an optimal policy. A hand-designed reward is necessary for each problem, which usually involves deep prior knowledge of complex physical systems. In addition, this approach requires numerical solvers during every iteration. These make such approaches inapplicable to real-world problems. 

We propose an alternative framework to solve PDE-constrained optimal control problems. Our method is divided into two phases: solution operator learning for PDE constraints (Phase 1) and searching for optimal control (Phase 2). During Phase 1, a neural network with a reconstruction structure is trained to approximate the PDE solutions. The optimal control problem can then be solved using the trained network with a reconstruction regularizer during Phase 2.
 The proposed method has the following four main contributions compared to the existing approaches in the literature:

\paragraph{Simple but effective regularizer}We introduce a novel regularizer to solve PDE-constrained optimal control problems effectively. We employ a reconstruction loss as a regularizer, which enables the control variable to converge correctly in Phase 2.

\paragraph{Application to various PDE-constrained control problems}We propose a non-problem-specific methodology for PDE-constrained optimal control problems. Our method is applied to problems with various types of PDEs: elliptic (Poisson, Stokes), hyperbolic (wave), non-linear parabolic (Burgers') equations. In addition, the control variable can be any type, such as the initial condition, the boundary condition, and a parameter in the governing equation.

\paragraph{Computational efficiency} Our surrogate model approximates a PDE solution with sufficiently high accuracy to search optimal controls while taking significantly less time for inference. Unlike the adjoint-based iterative methods, which require a heavy computation by PDE solvers for every iteration, our framework is useful when computation resources are limited, or a fast inference is required.

\paragraph{Flexibility} Our framework does not depend on the presence or absence of data. In the case that the full data pairs of control-to-state are available, the surrogate model for the PDE solution operator can be trained by a supervised loss during Phase 1. Meanwhile, if the physical systems are described in the form of PDEs, the residual norm of the PDE can be used to train the surrogate model without simulation data.

\section*{Related work}\label{realted}

\paragraph{Deep learning and PDEs}
There are two mainstream deep learning approaches to approximate solutions to the PDEs, i.e., using neural networks directly to parametrize the solution to the PDE and learning operators from the parameters of the PDEs to their solutions. A physics-informed neural network (PINN) was introduced in \cite{raissi2019physics}, which learns the neural network parameters to minimize the PDE residuals in the least-squares sense. In \cite{nabian2018deep}, \cite{son2021sobolev}, and \cite{weinan2018deep}, the authors suggested using a modified residual of the PDEs, and in \cite{han2018solving} and \cite{sirignano2018dgm}, the authors showed the possibility of solving high-dimensional PDEs. In \cite{hwang2020trend}, \cite{jo2020deep}, and \cite{sirignano2018dgm}, the authors prove a theorem on the approximation power of the neural network for an analytic solution to the PDEs. Next, we introduce another approach, operator learning, which is more closely related to our research.

\paragraph{Operator learning}
Operator learning using neural networks has been studied to approximate a PDE solution operator, which is nonlinear and complex in general. A universal approximation theorem for the operator was introduced in \cite{chen1993approximations}. Based on these results, in \cite{lu2019deeponet}, the authors developed a neural network called \textit{DeepONet}. In addition, mesh-based methods using convolutional neural networks (CNN) have been studied in many papers \cite{bhatnagar2019prediction, guo2016convolutional, khoo2017solving, zhu2018bayesian}. These studies used labeled data to train the operator networks. In \cite{bhatnagar2019prediction} and \cite{guo2016convolutional}, the authors used a CNN as a surrogate model of a computational fluid dynamics (CFD) solver. The authors showed that the surrogate models have a greater benefit in terms of speed than a CFD solver. In \cite{khoo2017solving} and \cite{zhu2018bayesian}, the authors developed the surrogate model for uncertainty quantification problems. Furthermore, the authors of \cite{zhu2019physics} proposed \textit{physics-constrained surrogate loss}, which can be calculated without labeled data. Li et al. proposed a resolution-invariant neural operator using a graph neural network \cite{li2020neural} and the fast Fourier transform \cite{li2020fourier}.

\paragraph{PDE-constrained control problem}
The most common approach in control theory is adjoint-based methods which give an efficient way to compute the gradient of forward maps with respect to the parameters \citep{borrvall2003topology,christofides2002nonlinear,lions1971optimal, mcnamara2004fluid, pironneau1974optimum,troltzsch2010optimal}. Several studies have suggested learning-based methods for control problems, such as \cite{de2018end}, \cite{hafner2019learning}, and \cite{NIPS2015_a1afc58c}. Regarding control problems associated with PDEs, the authors in \cite{holl2020learning} used a differentiable PDE solver to plan optimal trajectories and control fluid dynamics. They experimentally showed that their model enables long-term control with a fast inference time. Flow control problems, including vortex shedding and a drag reduction, were solved using reinforcement learning \cite{rabault2019artificial} or Koopman operator theory \cite{NEURIPS2018_2b0aa0d9}. One of the most interesting PDE-constrained optimization problems is an inverse problem, specifying unknown parameters in PDE systems given the observed data. There have recently been attempts to solve the problem by penalizing the parameter space or using a probabilistic approach \citep{jo2020deep, NEURIPS2019_28f0b864, ma2019probabilistic, pilozzi2018machine, ren2020benchmarking}. In particular, the approach in \cite{ren2020benchmarking} is similar to our study in that it learns the forward map first, but does not target the PDE problems. The authors in \cite{li2020fourier} showed that the PDE solution operator approximated by neural networks can be used in a Bayesian inverse problem.

\section*{Methodology}\label{method}

In this paper, we aim to solve PDE-constrained control problems. Let $M$ be a reflexive Banach space and $U$ and $V$ be Banach spaces. Formally, a PDE-constrained optimization problem can be written as follows:
\begin{equation}
    \min_{u \in U, m \in M} J(u, m) \quad \text{subject to } \quad F(u, m) = 0    
\end{equation}
where $J: U \times M \rightarrow \mathbb{R}$ is an objective function of interest, and $F: U \times M \rightarrow V$ is a system of PDEs, which governs the physics of the problem, possibly including initial and boundary conditions. Each space is a space of functions defined in a certain spatial or time domain. Here, $u$ is called a \textit{state} or PDE solution, and $m$ is a \textit{control} input. In the presence of control constraints, the problem is restricted to a set of admissible controls by $M_{ad} \subset M$, which is often assumed to be closed and convex. We remark that the control input can be configured in various forms, such as the values of a source term in a governing equation, or of the initial or boundary conditions. Our goal is to propose a general neural network based framework that is applicable to any type of PDEs and control inputs.

The optimal control of the Poisson equation can be considered a motivating example, a problem of specifying an unknown heat source to achieve a desired temperature profile. In this case, the control input $m$ indicates the values of the source term in the governing equation. The corresponding optimization problem is as follows:
\begin{equation}\label{poisson_control_objective}
  \min_{u \in H_0^1(\Omega), m \in L^2(\Omega)} \frac{1}{2}\int_{\Omega} (u-u_d)^2 \, dx  + \frac{\alpha}{2} \int_\Omega m^2 \, dx
\end{equation}
subject to the Poisson equation with zero Dirichlet boundary conditions
\begin{equation}\label{poisson_equation}
    \begin{cases}
    - \Delta u - m = 0 &\text{in}\ \Omega \\
    u = 0 &\text{on}\ \partial \Omega 
  \end{cases}    
\end{equation}
where $\Omega$ is the domain of interest, the state $u: \Omega \rightarrow \mathbb{R}$ is the unknown temperature, $u_d: \Omega \rightarrow \mathbb{R}$ is the given desired temperature, $\alpha$ is a penalty parameter, and $m: \Omega \rightarrow \mathbb{R}$ is the control function. Note that the penalty term should be distinguished from the regularization, which will be discussed in Section \ref{section3_2}. For practical purpose, $m$ is often restricted to $M_{ad}$, in which additional inequality constraints $m_a(x) \leq m(x) \leq m_b(x)$ are imposed. It is well-known that this problem is well-posed and has a unique solution \cite{hinze2008optimization, troltzsch2010optimal}.

We remark from the above example that the PDE solution $u$ can be thought as a function of $m$ with implicit relations $F(u, m) = 0$. In many cases, handling such complex, possibly nonlinear PDE constraints becomes the main difficulty when solving optimal control problems. One possible approach is to use surrogate models for PDE systems, approximating the explicit control-to-state mapping. For a concise notation, we denote the explicit solution expression by $u(m)$ and consider the reduced optimization problem
\begin{equation}
    \min_{m \in M_{ad}} \tilde{J}(m)
\end{equation}
with the reduced objective function $\tilde{J}(m) := J(u(m), m)$. This enables us to convert the constrained optimization problem into an unconstrained problem. We can then apply gradient-based optimization algorithms to obtain a locally optimal solution.

\begin{figure}[t]
\centering
\includegraphics[width=0.9\columnwidth]{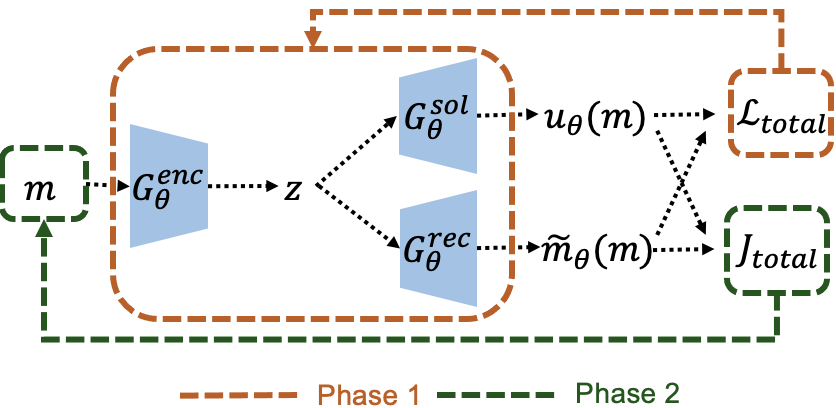} 
  \caption{Overview of our autoencoder model. During Phase 1, the parameter $\theta$ is updated, and during Phase 2, the control input $m$ is updated. }
  \label{architecture}
\end{figure}

In this study, we approximate the solution operators of PDE constraints through neural networks (Phase 1) and use them to search optimal controls for the given problems (Phase 2) through gradient descent. The two phases are described in detail in Section \ref{section3_1} and Section \ref{section3_2}, and summarized in Figure \ref{architecture}. In Section \ref{section3_3}, error estimates of the optimal controls are discussed under certain assumptions. Further, in Section \ref{section3_4}, we suggest modified architectures that are particularly efficient for time-dependent PDE constraints.

\subsection*{Phase 1: Solution operator learning for PDE constraints} \label{section3_1}
In the first step, a neural network is trained as a surrogate model for the PDE solution operator. We discretize the spatial domain into a uniform mesh to convert state $u$ and control input $m$ into image-like data. The surrogate model is then trained as an image-to-image regression. We suggest a variant autoencoder specialized for control problems. Our baseline model consists of a single encoder $G_\theta^{enc}$ for control input $m$ and two decoders $G_\theta^{sol}$, $G_\theta^{rec}$ corresponding to the state $u_\theta (m)=(G_\theta^{sol}\circ G_\theta^{enc})(m)$ and reconstruction $\tilde{m}_\theta (m)=(G_\theta^{rec}\circ G_\theta^{enc})(m)$ where $\theta$ is the set of all network parameters (Phase 1 in Figure \ref{architecture}). The reconstruction $\tilde{m}_\theta(m)$ plays an essential role in Phase 2. This will be described in Section \ref{section3_2}. Our method can be applied to the following two scenarios, \textit{data-driven} and \textit{data-free}.

\paragraph{Data-driven scenario}
%The loss function can be flexibly selected depending on the presence or absence of labeled data. 
In the case that the full data pairs of control-to-state are available, a supervised loss is a natural choice:
\begin{align}\label{loss_supervised}
  \mathcal{L}_{sup} = \frac{1}{N}\sum_{i=1}^N L(u_\theta (m_i), u_{i}),
\end{align}
where $\{(m_i, u_{i})\}_{i=1,\dots,N}$ is the observed data, and $L$ is a measure of the difference between two vectors. In our experiments, we used the $L^2$-relative error for $L$, namely $
    L(u, \tilde{u}) := \norm{u - \tilde{u}}_2 / \norm{\tilde{u}}_2
$

\paragraph{Data-free scenario} 
%However, 
In most real-world scenarios, complete data pairs of the control-to-state cannot be accessed because of expensive simulations. In these situations, one may utilize prior knowledge about the system of interest, which is often described in the form of PDEs. In that case, inspired by \cite{zhu2019physics}, the surrogate model can be trained by minimizing the residual norm of the PDE:
\begin{align}\label{loss_residual}
  \mathcal{L}_{res} = \frac{1}{N}\sum_{i=1}^N \norm{F(u_\theta (m_i), m_i)}^2,
\end{align}
where $\norm{\cdot}$ is the norm in the Banach space $V$. Here, we sampled the inputs $\{m_i\}_{i=1,\dots,N}$ in a set of admissible controls $M_{ad}$. This loss function imposes the physical law of the PDE constraint $F(u, m)=0$ to the surrogate model. For example, in the case of the Poisson equation \eqref{poisson_equation}, the residual norm of $F(u, m)$ can be expressed as
\begin{equation}
   \norm{F(u, m)}=\norm{- \Delta u - m}_{L^2(\Omega)}+\norm{u}_{L^2(\partial\Omega)}.
\end{equation}
When calculating residual $F(u_\theta (m_i), m_i)$ in $\mathcal{L}_{res}$, the spatial gradients can be approximated efficiently using a convolutional layer with a fixed kernel which consists of the finite difference coefficient, and the boundary condition can be enforced exactly (See Appendix \ref{appendix_poisson}).

Combining Eq. \eqref{loss_supervised} or Eq. \eqref{loss_residual} with the reconstruction loss, $\mathcal{L}_{rec} = \frac{1}{N}\sum_{i=1}^N L(\tilde{m}_\theta(m_i), m_i ),$ we used the total loss as $\mathcal{L}_{total} = \mathcal{L}_{sup} + \lambda_1 \mathcal{L}_{rec}$ or $\mathcal{L}_{res} + \lambda_1 \mathcal{L}_{rec}$ where $\lambda_1$ is a hyperparameter.

\subsection*{Phase 2: Searching for optimal control} \label{section3_2}
After the surrogate model is trained during Phase 1, the learned parameter $\theta^*$ is fixed. We cosider $m$ as a learnable parameter and denote the objective function as
\begin{align}
  J_{obj}(m) &:= J(u_{\theta^*}(m), m).
\end{align}
Because the surrogate model is differentiable, the gradient of the objective function with respect to the control input can be directly calculated. Then, $J_{obj}$ can be used as a loss function. If only $J_{obj}$ is minimized, one issue is that control input $m$ may converge to local optima outside the training domain of the surrogate model. This causes performance degradation of the surrogate model. To handle this situation, we employ the reconstruction loss as a regularizer, i.e., $J_{rec}(m) := L(\tilde{m}_{\theta^*} (m), m)$. $J_{rec}$ and $\mathcal{L}_{rec}$ in Phase 1 are similar, but different in that $\theta$ is updated during Phase 1 while $m$ is updated during Phase 2 with the fixed $\theta^*$.

\begin{figure}[t]
\centering
\includegraphics[width=0.9\columnwidth]{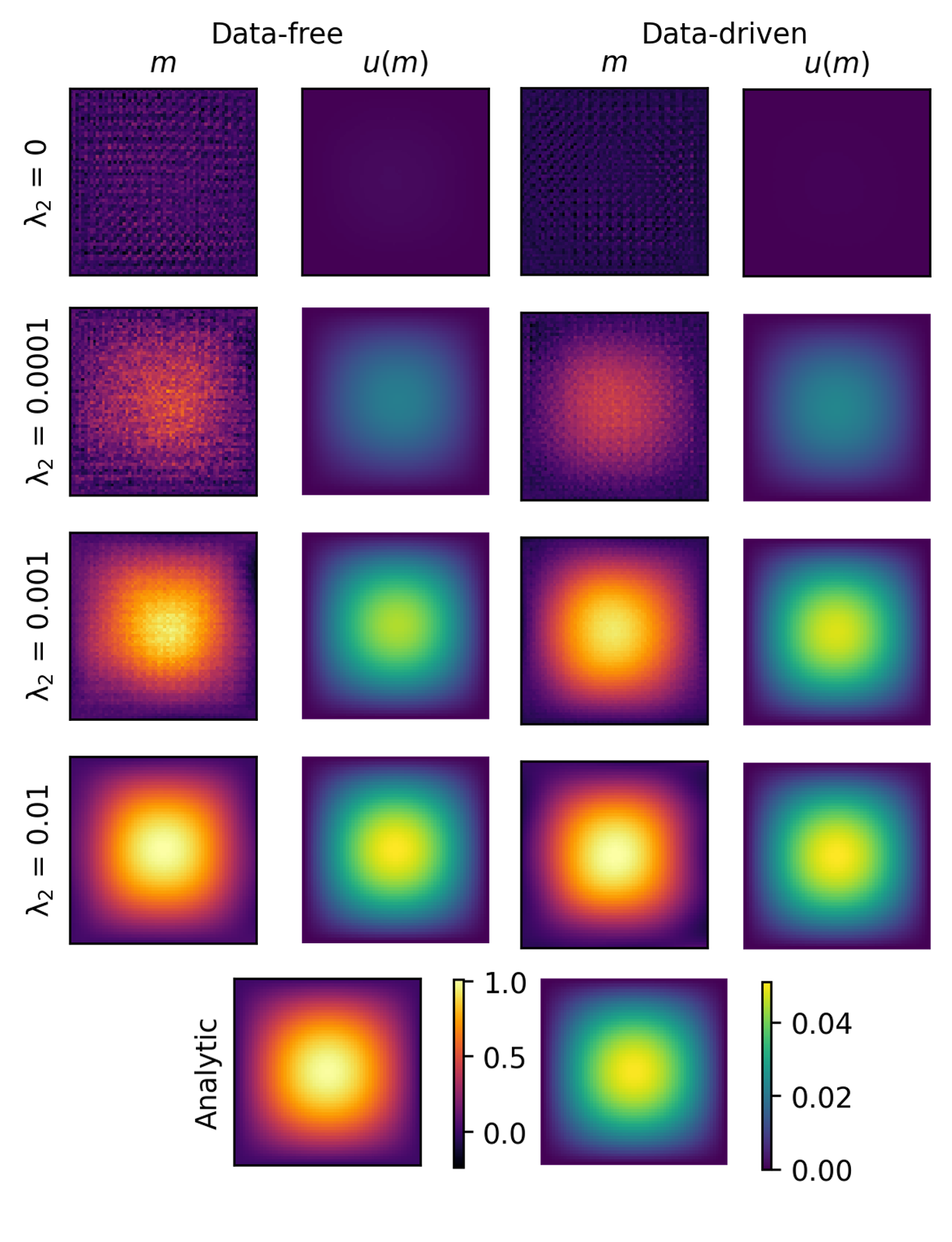} 
\caption{Control results according to different $\lambda_2$. First two columns are the results from data-free setting. The remaining columns are the results from data-driven setting. The last row shows the analytic optimal $m^*$ and $u^*$.}
\label{poisson_diff_lamb}
\end{figure}

Its regularizing effect can be interpreted in perspective of the variational autoencoder (VAE) \citep{kingma2013auto}. For this purpose, we consider the control input $m$ and latent variable $z$ as random vectors with prior density $p(z)$. A graphical model $p(m,z) = p(m|z)p(z)$ is given and induces an inequality given by 
\begin{multline}\label{vae_upperbound}
  -\log p(m) \\
   \leq \underbracket[0.8pt]{-\mathbb{E}_{q(z|m)}\left[\log p(m| z) \right] }_\text{\clap{reconstruction~}} +KL\left(q(z|m) |p(z) \right),
\end{multline} 
where $q(z|m)$ is an approximation of the posterior. $p(m)$ can be thought as the distribution of $m$ sampled during Phase 1 training. In this regard, we expect that minimizing the upper bound in Eq. \eqref{vae_upperbound} during Phase 2 keeps the likelihood $p(m)$ large enough, which implies that the updated control input $m$ keeps belonging to the training domain. In our experiments, we use a plain autoencoder, which models $q(z|\tilde{m})$ as a dirac distribution. In this case, the Kullback–Leibler divergence (KL) term is a constant with respect to $m$ and only the reconstruction term remains in the upper bound, which coincides with $J_{rec}(m)$ when $L$ is the $L^2$-loss and $p(m|z)$ is modeled through a Gaussian distribution. It implies that the role of the reconstruction regularizer in Phase 2 is to hold the control input $m$ in the region where the operator network works well. The experiments in Section \ref{experiment} show that the regularizer term greatly improves the performance of the optimal control learning, especially in Figure \ref{poisson_diff_lamb}.

Consequently, we set the following loss function to train the control problem:
\begin{equation}
  J_{total}(m)=J_{obj}(m)+\lambda_2 J_{rec}(m),
\end{equation}
where $\lambda_2$ is a hyperparameter.

\subsection*{Theoretical connection from Phase 1 to Phase 2} \label{section3_3}
Under some mild assumptions regarding a function space of neural network approximators, we derive the error estimates of the optimal control during Phase 2 in terms of the error that occurred during Phase 1. This provides the theoretical connection between the two separate phases. Although the following discussion is focused on our motivating example, i.e., a tracking-type problem with the Poisson equation, it can be adapted to other problems in a similar fashion.

In Eq. \eqref{poisson_equation}, we denote an exact solution operator by $S$, which satisfies $F(Sm, m) = 0$, and an approximate solution operator by $S_h$. In addition, let $\bmm$ be an exact optimal solution, and $\bmmh$ be the optimal solution inferred by our method. One of the main assumptions necessary to present our proposition is the Lipschitz continuity of the surrogate model $S_h$. Such an assumption has been well addressed in recent papers such as \cite{fazlyab2019efficient} amd \cite{NEURIPS2018_d54e99a6}. We then derive the following proposition:
\begin{prop}
  If $S_h: M \rightarrow U$ is Lipschitz continuous and approximated with error $\norm{S - S_h}_2 < \epsilon$, then the $L^2$ error for optimal control is estimated as
  $$
    \norm{\bmm - \bmmh}_2 < C \alpha^{-1} (1 + \alpha^{-1/2}) \norm{u_d}_2 \epsilon
  $$ 
  for a constant $C$.
\end{prop}
The observation $u_d$ and the penalty parameter $\alpha$ are the prescribed values. 
Meanwhile, $\norm{S - S_h}_2$, which is defined by
$
    \norm{S - S_h}_2 := \sup_{m\in L^2(\Omega),\ \norm{m} \leq 1} \norm{(S-S_h) m}_2,
$
can be thought of as a measurement of approximation and generalization of the operator learning. A sketch of proof is as follows: First we define the discrete version of the given optimization problem, which attains an optimal solution $m_h^*$. Then, by subtracting and modifying the first-order optimality condition for each problem, we can derive a proper 
upper bound for the error $\norm{\bmm - \bmmh}_2$. The detailed statements are given in Appendix \ref{appendix_theory}. 

The proposition provides the error estimates for the approximated optimal control input $\bmmh$ obtained using our method in Phase 2. 
If the surrogate model $S_h$ is trained to minimize the operator norm $\norm{S - S_h}_2$ within the given error tolerance $\epsilon$, then the error for optimal control $\norm{\bmm - \bmmh}_2$ can be estimated as in the proposition. This provides the connection between the two separate phases. 

\subsection*{Application to time-dependent PDEs} \label{section3_4}
We explain how our model can be extended and applied to time-dependent PDEs. We consider the time-dependent PDEs, in which the system can be written as
\begin{equation}
    \frac{\partial}{\partial t} u(t,\cdot)=F(u(t,\cdot),m(t,\cdot))
\end{equation}
for $t\in[0,T]$. We discretize the time-dependent PDE as follows:
\begin{equation}
    u^{t+\Delta t}=\mathcal{F}(u^t,m^t)
\end{equation}
where $u^t(\cdot):=u(t,\cdot)\in U$ and $m^t(\cdot):=m(t,\cdot)\in M$, in which $t=0,\Delta t, ...,(n-1)\Delta t$ with $T=n\Delta t$. Recently, the authors of \cite{lusch2018deep} employed a deep learning approach to discover representations of Koopman eigenfunctions from data. One of the key ideas in \cite{lusch2018deep} is that the time evolution of the eigenfunctions proceeds on the latent space between encoder and decoder. Our method can be extended for the time-dependent PDEs with inspiration from the idea in \cite{lusch2018deep}.

We use two autoencoders $H_\Theta$ and $G_\Theta$ for the state $u^t$ and control input $m^t$. The transition network $T_\Theta$ predicts the next time latent state $v^{t+\Delta t}$ from the latent state $v^t=H^{enc}(u^t)$ under the influence of the latent variable $g^t=G^{enc}(m^t)$. Therefore, the model propagates the state $u^t$ to the next time state $u^{t+\Delta t}$ under the influence of the control input $m^t$. The model can be used repeatedly $n$ times to generate a desired state $u^T$ from the given initial state $u^0$ and the given control inputs $m^0,\ m^{\Delta t},\ ...,\ m^{(n-1)\Delta t}$. The loss functions for Phases 1 and 2, and other details of the extended model are given in Appendix \ref{appendix_time_model}. 

\begin{figure*}[t]
\centering
\includegraphics[width=0.9\textwidth]{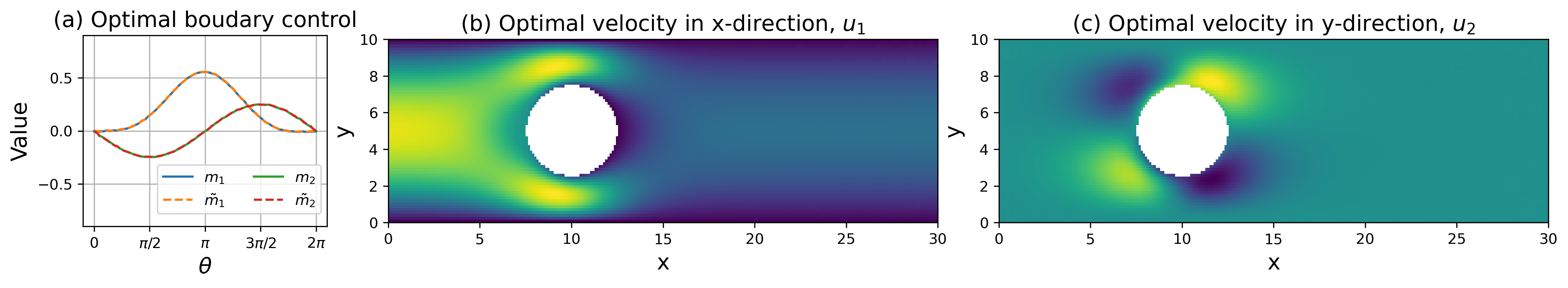} 
\caption{Illustration of the boundary control of the Stokes equation. (a) The control inputs $m_1$ and $m_2$ which corresond to the Dirichlet boundary values on the circle. The solid lines indicate the optimal control obtained by our method, and the dashed lines are the reconstructed control. (b), (c) The corresponding velocity fields $u_1$ and $u_2$ that minimize the drag energy.}
\label{stokes_control}
\end{figure*}

\section*{Experiment}\label{experiment}
In this section, we evaluate our method to handle a variety of PDE-constrained optimal control problems through four different examples. We first target the source control of the Poisson equation to illustrate the basic idea of our methodology. Two cases for operator learning, data-driven and data-free, will be considered to confirm the flexibility of our method. Next, we study the boundary control of the Stokes equation and the inverse design of the wave equation. Finally, the modified model architecture proposed in Section \ref{section3_4} will be verified using a nonlinear time-dependent PDE, Burgers' equation. For each experiment, we compare our method with the adjoint-based iterative method that is typically used for solving the optimal control problem effectively.

In this section, we focus on the results of optimal control during Phase 2. For the results of Phase 1, we reported the values of relative errors on test data in each experiment below. The visual results of the trained solution operator during Phase 1 are described in each subsection of Appendix \ref{appendix_experiment}. It shows that the surrogate models for the PDE solution operator in each control problem are well approximated with small relative errors, which is sufficient for use in Phase 2. Data generation and other details are given in Appendix \ref{appendix_experiment}.

\paragraph{Source control of the Poisson equation}
The Poisson equation is an elliptic PDE, also referred to a steady state heat equation. We consider the optimal control problem with control objective function \eqref{poisson_control_objective} subject to the Poisson equation under the Dirichlet boundary condition, which is described in equation \eqref{poisson_equation}. This is a fundamental PDE-constrained control problem. The problem aims to control the source term $m$ to make $u(m)$ similar to $u_d$ with $L^2$ penalization. During our experiment, $(x,y) \in \Omega = [0, 1]\times [0, 1]$, $\alpha = 10^{-6}$, and $u_d = \frac{1}{2 \pi^2} \sin \pi x \sin \pi y$. In this case, the analytic optimal control $m^*$ and the corresponding state $u^*$ are given by 
\begin{align}
  m^* = \frac{1}{1+ 4 \alpha \pi^4} \sin \pi x \sin \pi y, \quad
  u^* = \frac{1}{2 \pi^2} m^*.
\end{align}

Figure \ref{poisson_operator} in Appendix \ref{appendix_poisson} shows the training results for the solution operator during Phase 1. The solution operators are well approximated for both cases, data-driven (supervised loss, Eq. \eqref{loss_supervised}) and data-free (residual loss, Eq. \eqref{loss_residual}). The relative errors for test data are 0.0016 and 0.0080, respectively.

To verify the regularization effect of $J_{rec}$, we observe the change in optimal control obtained during Phase 2 depending on the regularizer coefficient $\lambda_2$. As shown in Figure \ref{poisson_diff_lamb}, with $\lambda_2=0$, i.e., no regularizer, the obtained optimal control is irregular and fails to converge to the analytic optimal $m^*$. As $\lambda_2$ increases, the optimal control becomes smoother and closer to $m^*$. This shows that the regularizer term greatly improves the performance of the optimal control learning. It makes the control input $m$ remain in the training domain where the operator network works well. In both cases, data-free and data-driven, the optimal control $m$ from our method is sufficiently close to analytic optimal control $m^*$ when $\lambda_2=0.01$. This phenomenon agrees with our expectation of the regularization effect discussed in Section \ref{section3_2}.

\paragraph{Boundary control of the stationary Stokes equation}
Consider the drag minimization problem of the two dimensional stationary Stokes equation:  
\begin{equation}\label{stokes_governing}
  \begin{cases}
        - \Delta u + \nabla p = 0 & \text{in } \Omega\\
    \text{div } u = 0 & \text{in } \Omega
  \end{cases}
\end{equation}\\
with boundary conditions
\begin{equation}
    \begin{cases}
        u = m & \text{on } \partial \Omega_{circle} \\
        u = f & \text{on } \partial \Omega_{in}, 
    \end{cases}
    \begin{cases}
        u = 0 & \text{on } \partial \Omega_{walls} \\ 
        p = 0 & \text{on } \partial \Omega_{out},
    \end{cases}
\end{equation}
where $\Omega$ is a rectangular domain with a circular obstacle inside (Each component is described in detail in Figure \ref{stokes_domain} in Appendix \ref{appendix_stokes}.). 
$u = [u_1, u_2]$ is the velocity, $p$ is the pressure, and $m = [m_1, m_2]$ is the control input, which corresponds to the Dirichlet boundary condition on the circle. The inflow boundary condition is given as $f(y) = y(10-y)/25$. This problem is interpreted as minimizing the drag from the flow by actively controlling the in/outflow on the circle boundary. 

\begin{figure}[t]
\centering
\includegraphics[width=0.9\columnwidth]{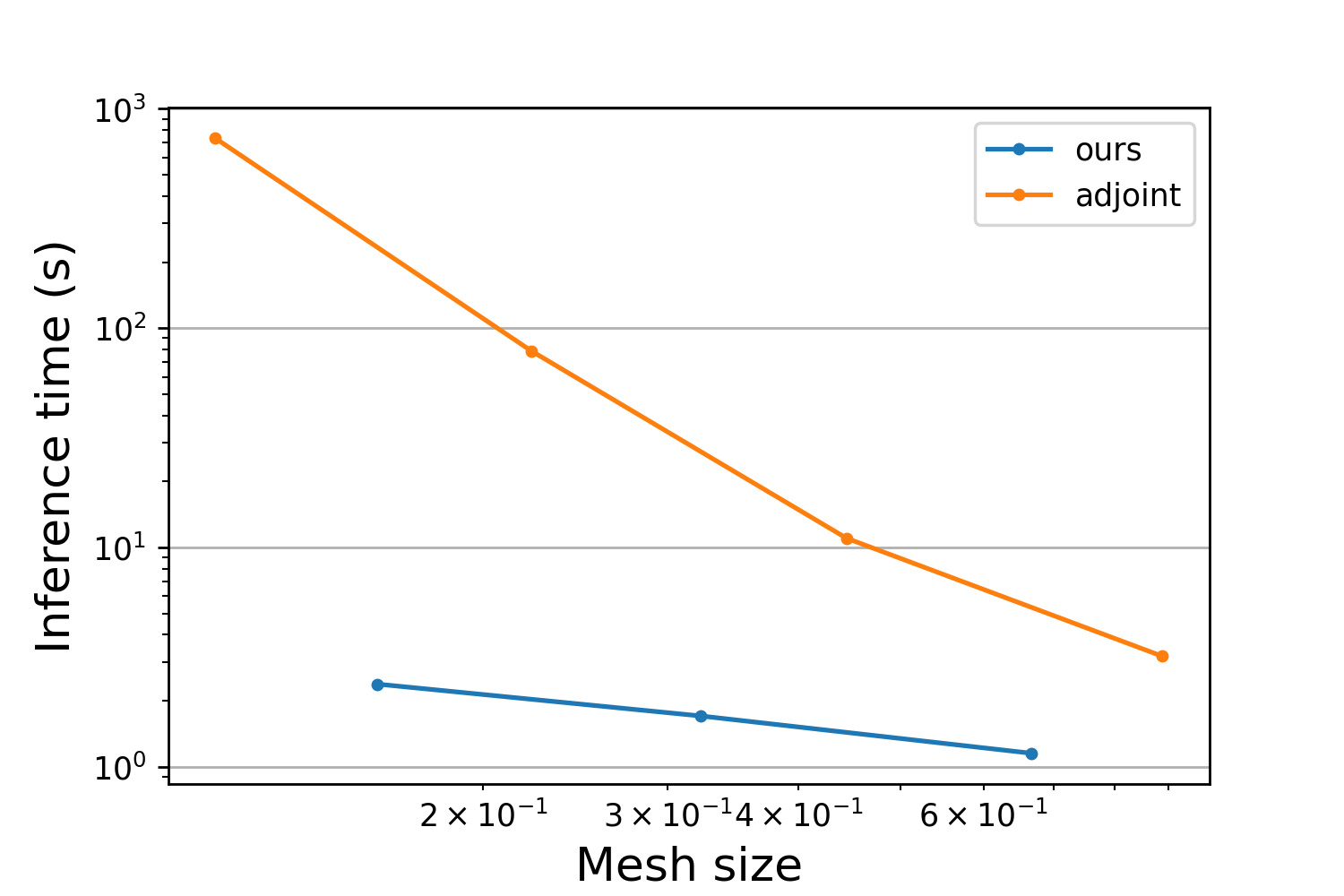} 
\caption{Comparison of the required inference time for the boundary control of the Stokes equation, using our method and the adjoint method. The results are plotted in the log-log scale. This shows that our method achieves better computational complexity than the adjoint method.}
\label{stokes_grid_vs}
\end{figure}

In Phase 1, the solution operator is well approximated with the relative error 0.0042 for test data. For the results of Phase 2, the left plot in Figure \ref{stokes_control} describes the obtained optimal control through our method. $m_1$ and $m_2$ are the $x$- and $y$-components of the optimal control, respectively, which are represented as functions of angle $\theta$ with respect to the center of the $\Omega_{circle}$. The right plot describes the corresponding velocity fields $[u_1, u_2]$ evaluated by the surrogate model. 
Figure \ref{stokes_grid_vs} summarizes the comparison of the inference time for optimal control (Phase 2) by varying the size of mesh when our method and the adjoint method are used. The mesh size in the $x$-axis means the maximum diameter of meshes used in each method. We remark that the complexity of the adjoint method increases much faster than our method as the mesh size increases. This is because the adjoint method requires heavy computation to obtain the exact gradient value of the objective function with respect to the control parameter.

\begin{figure*}[t]
\centering
\includegraphics[width=0.8\textwidth]{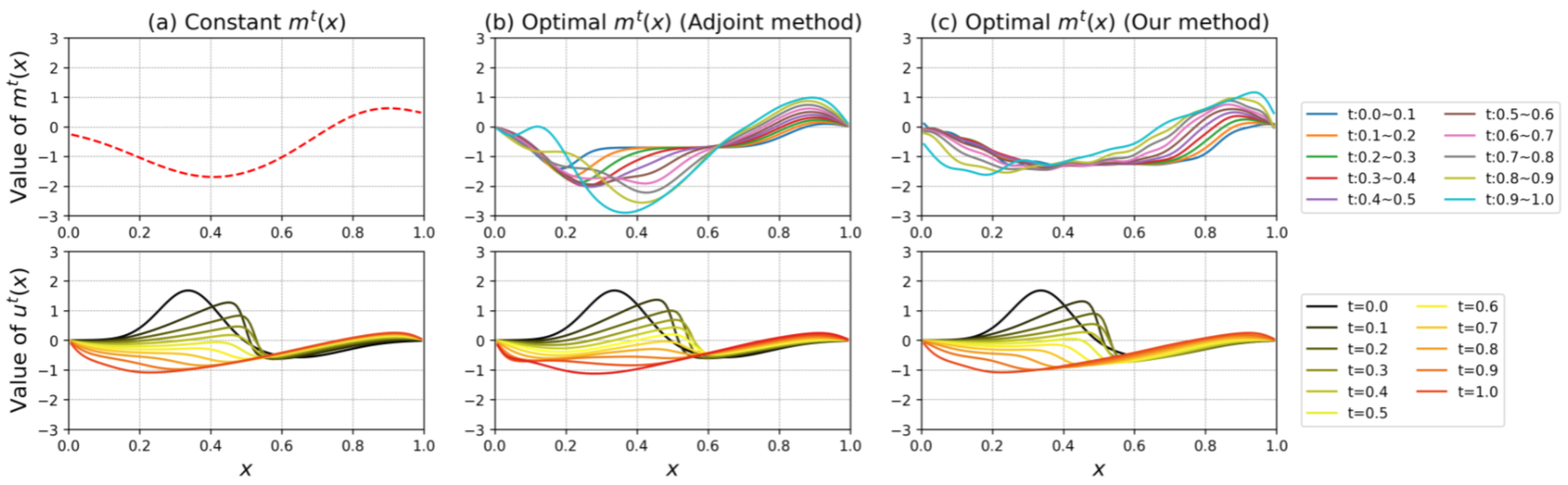} 
  \caption{The external force and trajectories with the external force using Burgers' equation. The three columns are (a) a constant external force, (b) an optimal external force using the adjoint method, (c) an optimal external force using our method. The first row shows the time-discretized optimal control inputs $m^t$ ($t=0,\Delta t, ...,(n-1)\Delta t$). The second row shows the time evolution of Burgers' equation with the external force in the first row.}
  \label{burgers_control_ex1}
\end{figure*}

\paragraph{Inverse design of nonlinear wave equation}
Given a target function $u_d(x)$ the inverse design problem aims to find the initial conditions that yield a solution $u(T, x) = u_d(x)$. The optimization problem is given by
\begin{equation}
  \min_{u, m} \frac{1}{2}\int_{\Omega} \left( u(T, x) - u_d(x) \right)^2 \, dx
\end{equation}
subject to the following nonlinear wave equation
\begin{equation}
    \begin{cases}
        \frac{\partial^2 u}{\partial t^2} - a^2 u_{xx} + f(u) = 0 &\text{in }[0, T]\times\Omega\\
        u(t, 0) = u(t, L) = 0 &\text{on } [0, T]\times\partial\Omega\\
        u(0, x) = 0, u_t(0, x) = m(x)  &\text{in }\{t=0\}\times\Omega
    \end{cases}
\end{equation}
where $\Omega = [0, L]$, and nonlinear source term $f(u) = u + u^3$. Here the initial condition $m(x)$ corresponds to control input. We choose $L=1$, $a=1/3$, and $T=5$. If the wave equation has a linear source term, the problem can be easily solved backward in time (time-reversibility). In our problem, however, such an approach fails owing to time-irreversible property caused by nonlinearity $f(u)$. Our method uses the surrogate model for direct mapping from initial to target state. 

In Phase 1, the solution operator is well approximated with the relative error 0.0065 for test data. For the results of Phase 2, Table \ref{wave_table} summarizes the results of the optimal control when using our method and the adjoint method. The objective function values of the two methods are comparable, but in terms of the computation time, our model significantly outperforms the adjoint method. This is because the surrogate model can infer fast the solution at $t=T$ for different control inputs, whereas the adjoint method needs to compute the forward and backward computations for each time step to reach the target time. In this regard, our method can be applied robustly even when the target time becomes longer.

\begin{table}[t]
\centering
\begin{tabular}{lll}
    \Xhline{1pt}
    
            & Objective     & Time (s)  \\ 
    \hline
    Ours    & 0.014 $\pm$ 0.005   & 0.210 $\pm$ 0.017  \\ 
    Adjoint & 0.012 $\pm$ 0.005   & 473.909 $\pm$ 43.622  \\ 
    \Xhline{1pt}
    
\end{tabular}
\caption{Optimal control results of the wave equation, repeated for 50 samples.}
\label{wave_table}
\end{table}

\begin{table}[t]
\centering
\begin{tabular}{lll}
    \Xhline{1pt}
    
            & Objective     & Time (s)  \\ 
    \hline
    Ours & 0.002 $\pm$ 0.002  &  7.714 $\pm$ 1.600   \\ 
    Adjoint    & 0.004 $\pm$ 0.003      & 12.965 $\pm$ 4.306      \\ 
    \Xhline{1pt}
    
\end{tabular}
\caption{Optimal control results of Burgers' equation, repeated for 100 samples.}
    \label{burgers_table}
\end{table}

\paragraph{Force control of Burgers' equation}
We use the extended model explained in Section \ref{section3_4} for the control problem to Burgers' equation. It describes the interaction between the effects of nonlinear convection and diffusion. Burgers' equation leads the shock wave phenomenon when the viscosity parameter has a small value. With the Dirichlet boundary condition, the 1D Burgers' equation with external force $m(t, x)$ reads as
\begin{equation}\label{burgers}
  \begin{cases}
    \frac{\partial u}{\partial t} = -u \cdot \frac{\partial u}{\partial x} + \nu \frac{\partial^2 u}{\partial x^2} + m(t, x) &\text{in }[0, T]\times\Omega\\
    u(t, x) = 0 \quad &\text{on } [0, T]\times\partial\Omega\\
    u(0, x) = u^0(t)
  \end{cases} 
\end{equation}
where $\nu$ is a viscosity parameter, and $u_0(t)$ is an initial condition. During the experiment, we consider the control problem minimizing
\begin{multline}
  \min_{u, m} \frac{1}{2}\int_{\Omega} |u(T, x) - u_d(x)|^2 \, dx \\
  + \frac{\alpha}{2} \int_{\Omega\times[0, T]} |m(t, x)|^2 \, dx dt
\end{multline}
subject to Burgers' equation \eqref{burgers}, given a target $u_d(x)$. In this case, we control an external force $m^t$. We set $\alpha=0.01$, $\Omega=[0,1]$, $T=1$, and $\Delta t=0.1$. Also, we set the viscosity parameter $\nu=0.01$ to generate a shock wave.

The solution operator is well trained during Phase 2 with a relative error 0.0020 for test data. Figure \ref{burgers_control_ex1} and Table \ref{burgers_table} show the result of our method compared to the adjoint method for optimal control (Phase 2). A time step size for the adjoint method is chosen as $\Delta t=0.01$ since the method does not converge when the time step size is set to the coarse grid ($\Delta t=0.1$) under our setting. In Figure \ref{burgers_control_ex1}, the trajectory in the second and third columns are scattered with less force than the constant external force in the first column. A quantitative comparison of our method to the adjoint method is shown in Table \ref{burgers_table}. Our framework takes less time compared to the adjoint method, while the objective function values are similar. Our surrogate model can mimic the time evolution of Burgers' equation in the coarse time grid. It makes our method infer faster than the adjoint method for the control optimization problem of time-dependent PDEs.

\section*{Conclusion}\label{conclusion}
We presented a general framework for solving PDE-constrained optimal control problems. We designed the surrogate models for PDE solution operators with a reconstruction structure. It allowed our model to solve the optimal control efficiently by adopting the reconstruction loss as a regularizer. The experimental results demonstrated that the proposed method has a significant gain in time complexity compared to the adjoint method. Also, our framework can be applied flexibly for both data-driven and data-free control problems.

Although our proposed method can achieve many benefits, it cannot completely replace the existing numerical methods. We believe that the two approaches can be complementary. The numerical method has an advantage in terms of accuracy, and our method is computationally efficient. In general, numerical methods slow down exponentially as the number of dimensions increases. We believe that a deep learning method can alleviate this issue.

%\clearpage

\section*{Acknowledgments}

This work was supported by the National Research Foundation of Korea (NRF) grant funded by the Korea government (MSIT) (NRF-2017R1E1A1A03070105, NRF-2019R1A5A1028324) and by Institute for Information \& Communications Technology Promotion (IITP) grant funded by the Korea government(MSIP) (No.2019-0-01906, Artificial Intelligence Graduate School Program (POSTECH)).

\bibliography{aaai22}

\begin{thebibliography}{47}
\providecommand{\natexlab}[1]{#1}

\bibitem[{Aln{\ae}s et~al.(2015)Aln{\ae}s, Blechta, Hake, Johansson, Kehlet,
  Logg, Richardson, Ring, Rognes, and Wells}]{alnaes2015fenics}
Aln{\ae}s, M.; Blechta, J.; Hake, J.; Johansson, A.; Kehlet, B.; Logg, A.;
  Richardson, C.; Ring, J.; Rognes, M.~E.; and Wells, G.~N. 2015.
\newblock The FEniCS project version 1.5.
\newblock \emph{Archive of Numerical Software}, 3(100).

\bibitem[{Berg and Nystr{\"o}m(2018)}]{berg2018unified}
Berg, J.; and Nystr{\"o}m, K. 2018.
\newblock A unified deep artificial neural network approach to partial
  differential equations in complex geometries.
\newblock \emph{Neurocomputing}, 317: 28--41.

\bibitem[{Bhatnagar et~al.(2019)Bhatnagar, Afshar, Pan, Duraisamy, and
  Kaushik}]{bhatnagar2019prediction}
Bhatnagar, S.; Afshar, Y.; Pan, S.; Duraisamy, K.; and Kaushik, S. 2019.
\newblock Prediction of aerodynamic flow fields using convolutional neural
  networks.
\newblock \emph{Computational Mechanics}, 64(2): 525--545.

\bibitem[{Borrvall and Petersson(2003)}]{borrvall2003topology}
Borrvall, T.; and Petersson, J. 2003.
\newblock Topology optimization of fluids in Stokes flow.
\newblock \emph{International journal for numerical methods in fluids}, 41(1):
  77--107.

\bibitem[{Bouchouev and Isakov(1999)}]{bouchouev1999uniqueness}
Bouchouev, I.; and Isakov, V. 1999.
\newblock Uniqueness, stability and numerical methods for the inverse problem
  that arises in financial markets.
\newblock \emph{Inverse problems}, 15(3): R95.

\bibitem[{Chen and Chen(1993)}]{chen1993approximations}
Chen, T.; and Chen, H. 1993.
\newblock Approximations of continuous functionals by neural networks with
  application to dynamic systems.
\newblock \emph{IEEE Transactions on Neural Networks}, 4(6): 910--918.

\bibitem[{Christofides and Chow(2002)}]{christofides2002nonlinear}
Christofides, P.~D.; and Chow, J. 2002.
\newblock Nonlinear and robust control of PDE systems: Methods and applications
  to transport-reaction processes.
\newblock \emph{Appl. Mech. Rev.}, 55(2): B29--B30.

\bibitem[{de~Avila Belbute-Peres et~al.(2018)de~Avila Belbute-Peres, Smith,
  Allen, Tenenbaum, and Kolter}]{de2018end}
de~Avila Belbute-Peres, F.; Smith, K.; Allen, K.; Tenenbaum, J.; and Kolter,
  J.~Z. 2018.
\newblock End-to-end differentiable physics for learning and control.
\newblock \emph{Advances in neural information processing systems}, 31:
  7178--7189.

\bibitem[{Egger and Engl(2005)}]{egger2005tikhonov}
Egger, H.; and Engl, H.~W. 2005.
\newblock Tikhonov regularization applied to the inverse problem of option
  pricing: convergence analysis and rates.
\newblock \emph{Inverse Problems}, 21(3): 1027.

\bibitem[{Fazlyab et~al.(2019)Fazlyab, Robey, Hassani, Morari, and
  Pappas}]{fazlyab2019efficient}
Fazlyab, M.; Robey, A.; Hassani, H.; Morari, M.; and Pappas, G.~J. 2019.
\newblock Efficient and accurate estimation of lipschitz constants for deep
  neural networks.
\newblock \emph{arXiv preprint arXiv:1906.04893}.

\bibitem[{Gunzburger(2002)}]{gunzburger2002perspectives}
Gunzburger, M.~D. 2002.
\newblock \emph{Perspectives in flow control and optimization}.
\newblock SIAM.

\bibitem[{Guo, Li, and Iorio(2016)}]{guo2016convolutional}
Guo, X.; Li, W.; and Iorio, F. 2016.
\newblock Convolutional neural networks for steady flow approximation.
\newblock In \emph{Proceedings of the 22nd ACM SIGKDD international conference
  on knowledge discovery and data mining}, 481--490.

\bibitem[{Hafner et~al.(2019)Hafner, Lillicrap, Fischer, Villegas, Ha, Lee, and
  Davidson}]{hafner2019learning}
Hafner, D.; Lillicrap, T.; Fischer, I.; Villegas, R.; Ha, D.; Lee, H.; and
  Davidson, J. 2019.
\newblock Learning latent dynamics for planning from pixels.
\newblock In \emph{International Conference on Machine Learning}, 2555--2565.
  PMLR.

\bibitem[{Han, Jentzen, and Weinan(2018)}]{han2018solving}
Han, J.; Jentzen, A.; and Weinan, E. 2018.
\newblock Solving high-dimensional partial differential equations using deep
  learning.
\newblock \emph{Proceedings of the National Academy of Sciences}, 115(34):
  8505--8510.

\bibitem[{Haslinger and M{\"a}kinen(2003)}]{haslinger2003introduction}
Haslinger, J.; and M{\"a}kinen, R.~A. 2003.
\newblock \emph{Introduction to shape optimization: theory, approximation, and
  computation}.
\newblock SIAM.

\bibitem[{Hinze et~al.(2008)Hinze, Pinnau, Ulbrich, and
  Ulbrich}]{hinze2008optimization}
Hinze, M.; Pinnau, R.; Ulbrich, M.; and Ulbrich, S. 2008.
\newblock \emph{Optimization with PDE constraints}, volume~23.
\newblock Springer Science \& Business Media.

\bibitem[{Holl, Koltun, and Thuerey(2020)}]{holl2020learning}
Holl, P.; Koltun, V.; and Thuerey, N. 2020.
\newblock Learning to control pdes with differentiable physics.
\newblock \emph{arXiv preprint arXiv:2001.07457}.

\bibitem[{Hwang et~al.(2020)Hwang, Jang, Jo, and Lee}]{hwang2020trend}
Hwang, H.~J.; Jang, J.~W.; Jo, H.; and Lee, J.~Y. 2020.
\newblock Trend to equilibrium for the kinetic Fokker-Planck equation via the
  neural network approach.
\newblock \emph{Journal of Computational Physics}, 419: 109665.

\bibitem[{Jo et~al.(2020)Jo, Son, Hwang, and Kim}]{jo2020deep}
Jo, H.; Son, H.; Hwang, H.~J.; and Kim, E.~H. 2020.
\newblock Deep neural network approach to forward-inverse problems.
\newblock \emph{Networks \& Heterogeneous Media}, 15(2): 247.

\bibitem[{Khoo, Lu, and Ying(2017)}]{khoo2017solving}
Khoo, Y.; Lu, J.; and Ying, L. 2017.
\newblock Solving parametric PDE problems with artificial neural networks.
\newblock \emph{arXiv preprint arXiv:1707.03351}.

\bibitem[{Kingma and Welling(2013)}]{kingma2013auto}
Kingma, D.~P.; and Welling, M. 2013.
\newblock Auto-encoding variational bayes.
\newblock \emph{arXiv preprint arXiv:1312.6114}.

\bibitem[{Li et~al.(2020{\natexlab{a}})Li, Kovachki, Azizzadenesheli, Liu,
  Bhattacharya, Stuart, and Anandkumar}]{li2020fourier}
Li, Z.; Kovachki, N.; Azizzadenesheli, K.; Liu, B.; Bhattacharya, K.; Stuart,
  A.; and Anandkumar, A. 2020{\natexlab{a}}.
\newblock Fourier neural operator for parametric partial differential
  equations.
\newblock \emph{arXiv preprint arXiv:2010.08895}.

\bibitem[{Li et~al.(2020{\natexlab{b}})Li, Kovachki, Azizzadenesheli, Liu,
  Bhattacharya, Stuart, and Anandkumar}]{li2020neural}
Li, Z.; Kovachki, N.; Azizzadenesheli, K.; Liu, B.; Bhattacharya, K.; Stuart,
  A.; and Anandkumar, A. 2020{\natexlab{b}}.
\newblock Neural operator: Graph kernel network for partial differential
  equations.
\newblock \emph{arXiv preprint arXiv:2003.03485}.

\bibitem[{Liang, Lin, and Koltun(2019)}]{NEURIPS2019_28f0b864}
Liang, J.; Lin, M.; and Koltun, V. 2019.
\newblock Differentiable Cloth Simulation for Inverse Problems.
\newblock In Wallach, H.; Larochelle, H.; Beygelzimer, A.; d\textquotesingle
  Alch\'{e}-Buc, F.; Fox, E.; and Garnett, R., eds., \emph{Advances in Neural
  Information Processing Systems}, volume~32. Curran Associates, Inc.

\bibitem[{Lions(1971)}]{lions1971optimal}
Lions, J. 1971.
\newblock \emph{Optimal Control of Systems Governed by Partial Differential
  Equations:}.
\newblock Grundlehren der mathematischen Wissenschaften in Einzeldarstellungen
  mit besonderer Ber{\"u}cksichtigung der Anwendungsgebiete. Springer-Verlag.
\newblock ISBN 9783540051152.

\bibitem[{Lu, Jin, and Karniadakis(2019)}]{lu2019deeponet}
Lu, L.; Jin, P.; and Karniadakis, G.~E. 2019.
\newblock Deeponet: Learning nonlinear operators for identifying differential
  equations based on the universal approximation theorem of operators.
\newblock \emph{arXiv preprint arXiv:1910.03193}.

\bibitem[{Lusch, Kutz, and Brunton(2018)}]{lusch2018deep}
Lusch, B.; Kutz, J.~N.; and Brunton, S.~L. 2018.
\newblock Deep learning for universal linear embeddings of nonlinear dynamics.
\newblock \emph{Nature communications}, 9(1): 1--10.

\bibitem[{Ma et~al.(2019)Ma, Cheng, Xu, Wen, and Liu}]{ma2019probabilistic}
Ma, W.; Cheng, F.; Xu, Y.; Wen, Q.; and Liu, Y. 2019.
\newblock Probabilistic representation and inverse design of metamaterials
  based on a deep generative model with semi-supervised learning strategy.
\newblock \emph{Advanced Materials}, 31(35): 1901111.

\bibitem[{McNamara et~al.(2004)McNamara, Treuille, Popovi{\'c}, and
  Stam}]{mcnamara2004fluid}
McNamara, A.; Treuille, A.; Popovi{\'c}, Z.; and Stam, J. 2004.
\newblock Fluid control using the adjoint method.
\newblock \emph{ACM Transactions On Graphics (TOG)}, 23(3): 449--456.

\bibitem[{Mitusch, Funke, and Dokken(2019)}]{mitusch2019dolfin}
Mitusch, S.~K.; Funke, S.~W.; and Dokken, J.~S. 2019.
\newblock dolfin-adjoint 2018.1: automated adjoints for FEniCS and Firedrake.
\newblock \emph{Journal of Open Source Software}, 4(38): 1292.

\bibitem[{Morton et~al.(2018)Morton, Jameson, Kochenderfer, and
  Witherden}]{NEURIPS2018_2b0aa0d9}
Morton, J.; Jameson, A.; Kochenderfer, M.~J.; and Witherden, F. 2018.
\newblock Deep Dynamical Modeling and Control of Unsteady Fluid Flows.
\newblock In Bengio, S.; Wallach, H.; Larochelle, H.; Grauman, K.;
  Cesa-Bianchi, N.; and Garnett, R., eds., \emph{Advances in Neural Information
  Processing Systems}, volume~31. Curran Associates, Inc.

\bibitem[{Nabian and Meidani(2018)}]{nabian2018deep}
Nabian, M.~A.; and Meidani, H. 2018.
\newblock A deep neural network surrogate for high-dimensional random partial
  differential equations.
\newblock \emph{arXiv preprint arXiv:1806.02957}.

\bibitem[{Paszke et~al.(2019)Paszke, Gross, Massa, Lerer, Bradbury, Chanan,
  Killeen, Lin, Gimelshein, Antiga et~al.}]{paszke2019pytorch}
Paszke, A.; Gross, S.; Massa, F.; Lerer, A.; Bradbury, J.; Chanan, G.; Killeen,
  T.; Lin, Z.; Gimelshein, N.; Antiga, L.; et~al. 2019.
\newblock Pytorch: An imperative style, high-performance deep learning library.
\newblock \emph{arXiv preprint arXiv:1912.01703}.

\bibitem[{Pilozzi et~al.(2018)Pilozzi, Farrelly, Marcucci, and
  Conti}]{pilozzi2018machine}
Pilozzi, L.; Farrelly, F.~A.; Marcucci, G.; and Conti, C. 2018.
\newblock Machine learning inverse problem for topological photonics.
\newblock \emph{Communications Physics}, 1(1): 1--7.

\bibitem[{Pironneau(1974)}]{pironneau1974optimum}
Pironneau, O. 1974.
\newblock On optimum design in fluid mechanics.
\newblock \emph{Journal of Fluid Mechanics}, 64(1): 97--110.

\bibitem[{Rabault et~al.(2019)Rabault, Kuchta, Jensen, R{\'e}glade, and
  Cerardi}]{rabault2019artificial}
Rabault, J.; Kuchta, M.; Jensen, A.; R{\'e}glade, U.; and Cerardi, N. 2019.
\newblock Artificial neural networks trained through deep reinforcement
  learning discover control strategies for active flow control.
\newblock \emph{Journal of fluid mechanics}, 865: 281--302.

\bibitem[{Raissi, Perdikaris, and Karniadakis(2019)}]{raissi2019physics}
Raissi, M.; Perdikaris, P.; and Karniadakis, G.~E. 2019.
\newblock Physics-informed neural networks: A deep learning framework for
  solving forward and inverse problems involving nonlinear partial differential
  equations.
\newblock \emph{Journal of Computational Physics}, 378: 686--707.

\bibitem[{Ren, Padilla, and Malof(2020)}]{ren2020benchmarking}
Ren, S.; Padilla, W.; and Malof, J. 2020.
\newblock Benchmarking deep inverse models over time, and the neural-adjoint
  method.
\newblock \emph{arXiv preprint arXiv:2009.12919}.

\bibitem[{Sirignano and Spiliopoulos(2018)}]{sirignano2018dgm}
Sirignano, J.; and Spiliopoulos, K. 2018.
\newblock DGM: A deep learning algorithm for solving partial differential
  equations.
\newblock \emph{Journal of computational physics}, 375: 1339--1364.

\bibitem[{Sokolowski and Zol{\'e}sio(1992)}]{sokolowski1992introduction}
Sokolowski, J.; and Zol{\'e}sio, J.-P. 1992.
\newblock Introduction to shape optimization.
\newblock In \emph{Introduction to Shape Optimization}, 5--12. Springer.

\bibitem[{Son et~al.(2021)Son, Jang, Han, and Hwang}]{son2021sobolev}
Son, H.; Jang, J.~W.; Han, W.~J.; and Hwang, H.~J. 2021.
\newblock Sobolev Training for the Neural Network Solutions of PDEs.
\newblock \emph{arXiv preprint arXiv:2101.08932}.

\bibitem[{Tr{\"o}ltzsch(2010)}]{troltzsch2010optimal}
Tr{\"o}ltzsch, F. 2010.
\newblock \emph{Optimal control of partial differential equations: theory,
  methods, and applications}, volume 112.
\newblock American Mathematical Soc.

\bibitem[{Virmaux and Scaman(2018)}]{NEURIPS2018_d54e99a6}
Virmaux, A.; and Scaman, K. 2018.
\newblock Lipschitz regularity of deep neural networks: analysis and efficient
  estimation.
\newblock In Bengio, S.; Wallach, H.; Larochelle, H.; Grauman, K.;
  Cesa-Bianchi, N.; and Garnett, R., eds., \emph{Advances in Neural Information
  Processing Systems}, volume~31. Curran Associates, Inc.

\bibitem[{Watter et~al.(2015)Watter, Springenberg, Boedecker, and
  Riedmiller}]{NIPS2015_a1afc58c}
Watter, M.; Springenberg, J.; Boedecker, J.; and Riedmiller, M. 2015.
\newblock Embed to Control: A Locally Linear Latent Dynamics Model for Control
  from Raw Images.
\newblock In Cortes, C.; Lawrence, N.; Lee, D.; Sugiyama, M.; and Garnett, R.,
  eds., \emph{Advances in Neural Information Processing Systems}, volume~28.
  Curran Associates, Inc.

\bibitem[{Weinan and Yu(2018)}]{weinan2018deep}
Weinan, E.; and Yu, B. 2018.
\newblock The deep Ritz method: a deep learning-based numerical algorithm for
  solving variational problems.
\newblock \emph{Communications in Mathematics and Statistics}, 6(1): 1--12.

\bibitem[{Zhu and Zabaras(2018)}]{zhu2018bayesian}
Zhu, Y.; and Zabaras, N. 2018.
\newblock Bayesian deep convolutional encoder--decoder networks for surrogate
  modeling and uncertainty quantification.
\newblock \emph{Journal of Computational Physics}, 366: 415--447.

\bibitem[{Zhu et~al.(2019)Zhu, Zabaras, Koutsourelakis, and
  Perdikaris}]{zhu2019physics}
Zhu, Y.; Zabaras, N.; Koutsourelakis, P.-S.; and Perdikaris, P. 2019.
\newblock Physics-constrained deep learning for high-dimensional surrogate
  modeling and uncertainty quantification without labeled data.
\newblock \emph{Journal of Computational Physics}, 394: 56--81.

\end{thebibliography}

\clearpage
\appendix
\section{Theoretical connection from Phase 1 to Phase 2}\label{appendix_theory}
Under some mild assumptions on a function space of neural network approximators, we derive the error estimates of optimal control using approximate solution operators: If the solution operator is well approximated, then the optimal control can be inferred with a small error. Although the discussion is focused on our motivating example, i.e., a tracking-type problem with the Poisson equation, it can be adapted to other problems. In this case, we have $U = H^1_0(\Omega)$, $V = H^{-1}(\Omega)$, $M=L^2(\Omega)$. The states and control inputs are considered as piece-wise constant functions on an equidistant grid of $\Omega$ with mesh size $h$. Let $M_h \subset M$  and $U_h \subset U$ be the sets of all functions that are constant in each cell. Denote an (analytic) solution operator by $S: M\rightarrow U$, which satisfies $F(Sm, m) = 0$, and an approximate solution operator by $S_h: M_h \rightarrow U_h$. We then assume that the surrogate model $S_h$ is Lipschitz continuous, i.e., there exists a constant $C$ satisfying
\begin{equation}
  \norm{S_h (m_1 - m_2)}_2 \leq C \norm{m_1 - m_2}_2
\end{equation}
for any control inputs $m_1, m_2 \in M_h$. For the unconstrained case, the optimal control problem above is rewritten as
\begin{equation}
  \min_{m \in L^2(\Omega)} \tilde{J}(m) = \frac{1}{2} \norm{Sm - u_d}_2^2 + \frac{\alpha}{2} \norm{m}_2^2. \tag{$P$}
\end{equation}
The optimal control problem using the approximated solution operator $S_h$ is written as
\begin{equation}
  \min_{m \in L^2(\Omega) \cap M_h} \tilde{J}_h(m) = \frac{1}{2} \norm{S_h m - u_d}_2^2 + \frac{\alpha}{2} \norm{m}_2^2. \tag{$P_h$}
\end{equation}
We assume that the above problems attain optimal solution, namely, $\bmm = \arg \min \tilde{J}(m)$ and $\bmmh = \arg \min \tilde{J}_h(m)$. Now we will prove the main theoretical property. 

\begin{proof}[Proof of Proposition in Section \ref{section3_3}]
  It is easy to derive the necessary conditions for both optimal controls $\bmm$ and $\bmmh$, respectively:
\begin{gather*}
    \tilde{J}'(\bmm) = S^\dagger (S\bmm - u_d) + \alpha \bmm = 0,\\
    \tilde{J}_h'(\bmmh) = S_h^\dagger (S_h\bmmh - u_d) + \alpha \bmmh = 0, 
\end{gather*}
where $S^\dagger$ and $S_h^\dagger$ denote the adjoint operators. 

The error $\norm{\bmm - \bmmh}_2$ can then be estimated as follows: First, multiplying $\bmm - \bmmh$ on both equations and subtracting the results, we have
\begin{multline}
  (S^\dagger (S\bmm - u_d) - S_h^\dagger (S_h\bmmh - u_d), \bmm - \bmmh)_2 \\
  + \alpha \norm{\bmm - \bmmh}_2^2 = 0.
\end{multline}
Here, $(\cdot, \cdot)_2$ denotes the inner product in $L^2(\Omega)$. We add and subtract $(S_h^\dagger S_h \bmm, \bmm - \bmmh)_2$ on the first term, yielding 

\begin{multline}
  \norm{S_h (\bmm - \bmmh)}_2^2 + \alpha \norm{\bmm - \bmmh}_2^2 \\
    = \underbracket[0.8pt]{(u_d, (S - S_h) (\bmm - \bmmh))_2}_\text{\clap{\romanone~}} \\
    + \underbracket[0.8pt]{((S_h^\dagger S_h - S^\dagger S)\bmm, \bmm - \bmmh)_2}_\text{\clap{\romantwo~}}. 
\end{multline}
It is clear that
\begin{equation}
  \left| \text{\romanone} \right| \leq \norm{u_d}_2 \norm{S - S_h}_2 \norm{\bmm - \bmmh}_2.
\end{equation}
Here, $\norm{S - S_h}_2$ is the natural operator norm defined by
\begin{equation} \label{opnorm}
  \norm{S - S_h}_2 := \sup_{\substack{m\in L^2(\Omega) \\ \norm{m} \leq 1}} \norm{(S-S_h) m}_2.
\end{equation}
For \romantwo, adding and subtracting $(S_h^\dagger S \bmm, \bmm - \bmmh)_2$, we have
\begin{multline}
    \left| \text{\romantwo} \right| \leq \left| (S_h^\dagger (S_h - S) \bmm, \bmm - \bmmh)_2\right| \\
    + \left| ((S_h^\dagger - S^\dagger) S \bmm, \bmm - \bmmh)_2\right| \nonumber \\
    \leq \norm{(S_h - S) \bmm}_2 \norm{S_h (\bmm - \bmmh)}_2 \\
    + \norm{S \bmm}_2 \norm{(S_h - S) (\bmm - \bmmh)}_2 \nonumber \\
    \leq C \norm{\bmm}_2 \norm{S_h - S}_2 \norm{\bmm - \bmmh}_2
\end{multline}
based on the continuity of operators (namely, $\norm{S \bmm}_2 \leq C \norm{\bmm}_2$ and $\norm{S_h (\bmm - \bmmh)}_2 \leq C \norm{\bmm - \bmmh}_2$). Combining two estimates, we have an $L^2$ error estimate
\begin{equation}
  \norm{\bmm - \bmmh}_2 \leq C \alpha^{-1} (\norm{u_d}_2 + \norm{\bmm}_2) \norm{S - S_h}_2. 
\end{equation}
Finally, by optimality, $\norm{\bmm}_2 \leq \alpha^{-1/2} \norm{u_d}_2$, and the last term is estimated as  
\begin{align}
  \norm{\bmm - \bmmh}_2 &\leq C \alpha^{-1} (1 + \alpha^{-1/2}) \norm{u_d}_2 \norm{S - S_h}_2 \nonumber \\
  &< C \alpha^{-1} (1 + \alpha^{-1/2}) \norm{u_d}_2 \epsilon.
\end{align}
\end{proof}
We remark that the operator norm $\norm{S - S_h}_2$ can be thought as a measurement of approximation and generalization of the operator learning. In Phase 1, the surrogate model is trained to minimize the operator norm $\norm{S - S_h}_2$. If the model is within the given error tolerance, $\norm{S - S_h}_2 < \epsilon$, the proposition guarantees that the optimal control can be inferred with less error than $C \alpha^{-1} (1 + \alpha^{-1/2}) \norm{u_d}_2 \epsilon$. This provides the connection between two separate phases.

\section{Details of extended model to time-dependent PDE}\label{appendix_time_model}

\begin{figure*}[t]
\centering
\includegraphics[width=0.8\textwidth]{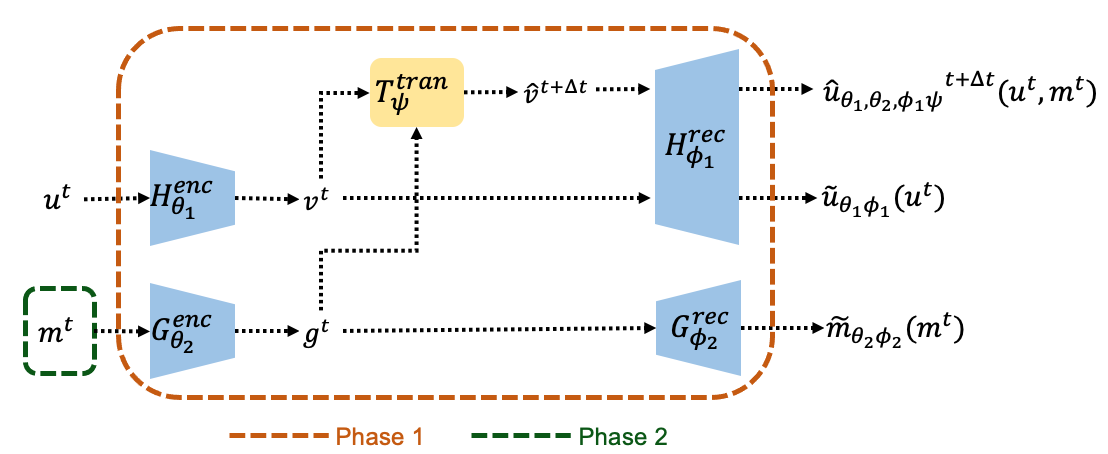} 
  \caption{Overview of our extended autoencoder model for time-dependent PDE. We denote the all network parameters as $\Theta:=[\theta_1,\theta_2,\phi_1,\phi_2,\psi]$. During Phase 1, the parameter $\Theta$ is updated and during Phase 2, the control input $m^t$ ($t=0,\Delta t, ...,(n-1)\Delta t$) is updated.}
  \label{architecture_time}
\end{figure*}

In this section, we explain in detail how our model can be extended and applied to the time-dependent PDE. We propose a model architecture that learns the dynamics on latent variables as a linear representation. The state $u^t$ at time $t$ is propagated to the next state $u^{t+\Delta t}$ in a linear framework. To this end, the two autoencoders, $H_{\theta_1,\phi_1}$ and $G_{\theta_2,\phi_2}$, are used for the state $u^t$ and the control input $m^t$ as in Figure \ref{architecture_time}. The encoder $H_{\theta_1}^{enc}$ takes the state $u^t$ as input to obtain the latent state $v^t(=H_{\theta_1}^{enc}(u^t))$. The decoder $H_{\phi_1}^{rec}$ is used to reconstruct the state $u^t$ from the latent state $v^t$. We denote the reconstruction of the state $u^t$ as $\tilde{u}_{\theta_1,\phi_1}(u^t)$.

For the control input $m^t$, we also use the another autoencoder structure $G_{\theta_2,\phi_2}$ to make $\tilde{m}_{\theta_2,\phi_2}(m^t)$, which plays a key role in Phase 2. The latent variable of the control input $m^t$ is denoted as $g^t(=G_{\theta_2}^{enc}(m^t))$. Then, the network $T_\psi^{tran}$ propagates the latent state $v^t$ to the next time $\hat{v}^{t+\Delta t}(=T_\psi^{tran}(v^t,g^t))$. The two latent images $v^t$ and $g^t$ are passed to the transition network $T_\psi^{tran}$ to get the latent state $\hat{v}^{t+\Delta t}$ at the next time $t+\Delta t$. The desired next state $\hat{u}^{t+\Delta t}(=H_{\phi_1}^{rec}(\hat{v}^{t+\Delta t}))$ is generated using the decoder $H_{\phi_1}^{rec}$. We also use another expression for the next state $\hat{u}^{t+\Delta t}$ as $\hat{u}_{\theta_1,\theta_2,\phi_1,\psi}^{t+\Delta t}(u^t,m^t)$.

Given observed data $\{ (u_i^{j\Delta t}, m_i^{j\Delta t}, u_i^{(j+1)\Delta t} )\}$$(1\leq i\leq N, \,0\leq j \leq n-1)$, we use the modified loss function
\begin{equation}\label{loss_burgers_phase1}
  \mathcal{L}_{total} = \mathcal{L}_{sup} + \lambda_1 \mathcal{L}_{rec}
\end{equation} for Phase 1, where
\begin{equation}
\mathcal{L}_{sup} := \frac{1}{N}\sum_{i=1}^N \sum_{j=1}^{n} L(\hat{v}_i^{j\Delta t}, v^{j\Delta t}_{i}) + L(\hat{u}_i^{j\Delta t}, u_i^{j\Delta t}),
\end{equation}
and
\begin{multline}
  \mathcal{L}_{rec} := \frac{1}{N}\sum_{i=1}^N \bigg(\sum_{j=0}^{n} L(\tilde{u}_{\theta_1, \phi_1}(u_i^{j\Delta t}), u_i^{j\Delta t})\\
  + \sum_{j=0}^{n-1}L(\tilde{m}_{\theta_2, \phi_2}(m_i^{j\Delta t}), m_i^{j\Delta t})\bigg).
\end{multline}

To search the optimal control input $m=\{m^{j\Delta t}\}_{j=0}^{n-1}$ during Phase 2, we use the objective function
\begin{equation}\label{loss_burgers_phase2}
  J_{total}(m)=J_{obj}(m)+\lambda_2 J_{rec}(m),
\end{equation}
with
\begin{equation}
J_{obj}(m) := J\left(\hat{u}(m), m\right)
\end{equation}
where $\hat{u}(m)=\{\hat{u}^{j\Delta t}(m)\}_{j=1}^{n}$ and
\begin{equation}\notag
    \hat{u}^{j\Delta t}(m)=\hat{u}_{\theta_1^*,\theta_2^*,\phi_1^*,\psi^*}^{j\Delta t}(\hat{u}^{(j-1)\Delta t},m^{(j-1)\Delta t}),
\end{equation} and with
\begin{multline}
J_{rec}(m) := \sum_{j=1}^{n} L(\tilde{u}_{\theta_1^*, \phi_1^*}(\hat{u}^{j\Delta t}(m)), \hat{u}^{j\Delta t}(m))\\
+ \sum_{j=0}^{n-1} L(\tilde{m}_{\theta_2^*, \phi_2^*}(m^{j\Delta t}), m^{j\Delta t}).
\end{multline}

We choose the relative error for the measure $L(\cdot,\cdot)$.

\begin{figure*}[t]
\centering
\includegraphics[width=0.7\textwidth]{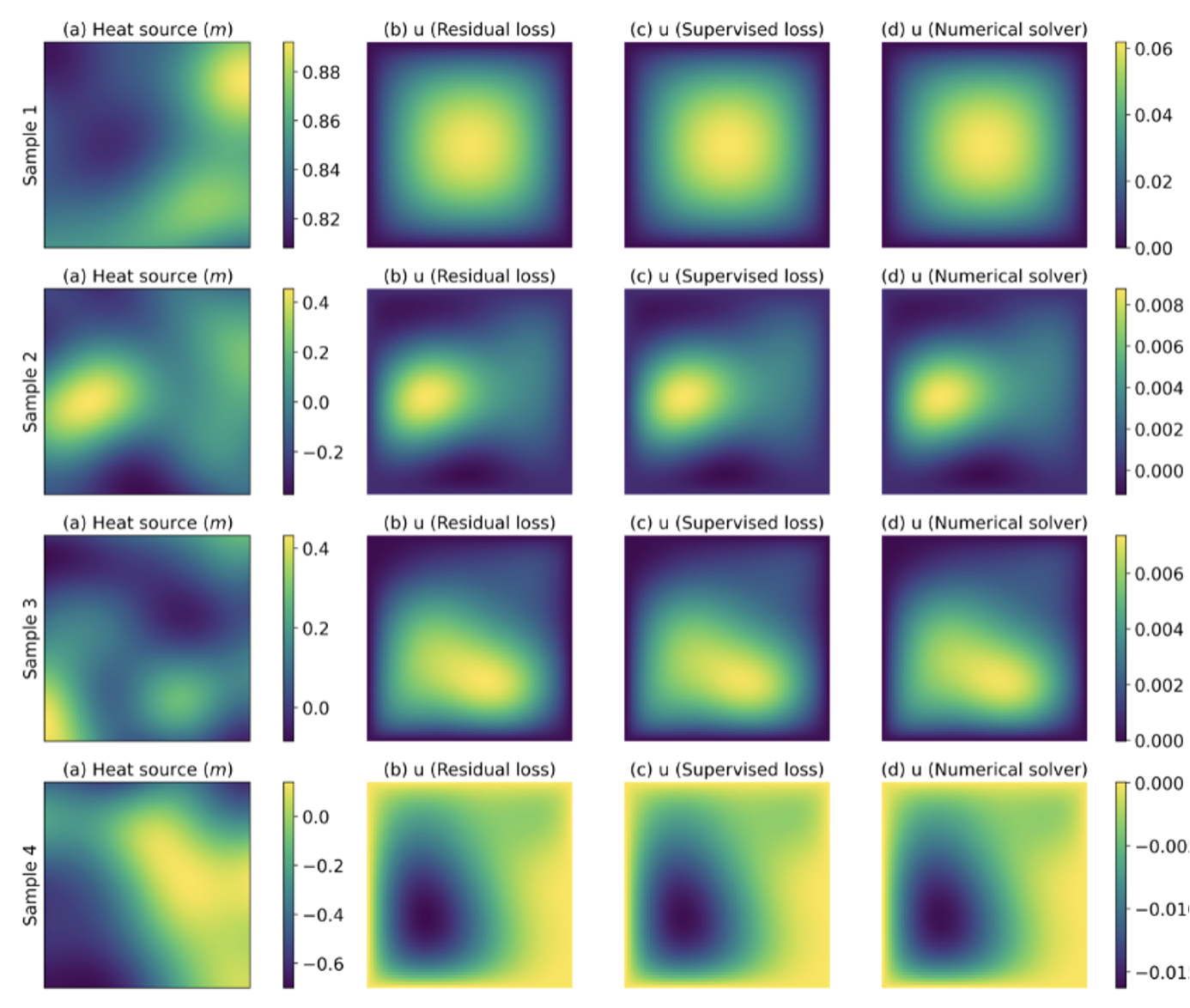} 
  \caption{Examples of input and output of operator network for the Poisson equation. (a) Randomly generated input heat source $m$. (b), (c) Predicted $u$ by the models trained with the residual loss and the supervised loss, repectively. (d) Predicted $u$ from numerical method.}
  \label{poisson_operator}
\end{figure*}

\begin{table*}[t]
\centering
\begin{tabular}{clrrrr}

  \Xhline{1pt}
Initial guess &
  \multicolumn{1}{c}{Method} &
  \multicolumn{1}{c}{Rel($m$, $m^*$)} &
  \multicolumn{1}{c}{Rel($u$, $u^*$)} &
  \multicolumn{1}{c}{Objective function} &
  \multicolumn{1}{c}{Time (s)} \\ \hline
\multirow{3}{*}{$m_{init}= xy(1-x)(1-y)$} & Residual loss    & $2.531\times10^{-2}$ & $6.098\times10^{-4}$ & $1.238\times10^{-7}$ & 1.176  \\
                                  & Supervised loss  & $2.131\times10^{-2}$ & $7.101\times10^{-4}$ & $1.253\times10^{-7}$ & 1.711 \\
                                  & Adjoint method & $\pmb{2.072\times10^{-2}}$ & $\pmb{5.021\times10^{-4}}$ & $\pmb{1.222\times10^{-7}}$ & \textbf{0.460}  \\ \hline
\multirow{3}{*}{$m_{init} = x + y$}      & Residual loss    & $\pmb{2.631\times10^{-2}}$ & $\pmb{5.070\times10^{-4}}$ & $\pmb{1.238\times10^{-7}}$ & \textbf{2.336}  \\
                                  & Supervised loss  & $3.320\times10^{-2}$ & $8.920\times10^{-4}$ & $1.253\times10^{-7}$ & 3.526 \\
                                  & Adjoint method & $1.615\times10^{-1}$ & $2.715\times10^{-3}$ & $1.275\times10^{-7}$ & 2.378  \\ \hline
\multirow{3}{*}{$m_{init} = xy$}         & Residual loss    & $\pmb{2.738\times10^{-2}}$ & $\pmb{5.667\times10^{-4}}$ & $\pmb{1.239\times10^{-7}}$ & \textbf{1.789}  \\
                                  & Supervised loss  & $4.721\times10^{-2}$ & $9.225\times10^{-4}$ & $1.246\times10^{-7}$ & 3.412 \\
                                  & Adjoint method & $1.291\times10^{-1}$ & $1.423\times10^{-3}$ & $1.246\times10^{-7}$ & 2.167 \\ 
  \Xhline{1pt}
\end{tabular}
\caption{Comparison of control results from different initial $m_{init}$. Rel$(\cdot,\cdot)$ denotes the relative error.}
\label{poisson_table}
\end{table*}

\section{Experimental details}\label{appendix_experiment}
In this section, we give further descriptions for the experiments of Section \ref{experiment}. 
We used the deep learning framework PyTorch \cite{paszke2019pytorch} (BSD license) for our neural network implementations. We compare our method and the adjoint-based iterative method. The adjoint method is implemented using the finite element framework FEniCS \cite{alnaes2015fenics} (LGPLv3 license) and its algorithmic differentiation tool dolfin-adjoint \cite{mitusch2019dolfin} (LGPLv3 license). We use AMD Ryzen 7 5800X processor and NVIDIA GeForce RTX 3080 10GB GPU. During Phase 1, the GPU is used for operator learning. On the other hand, during Phase 2, the CPU is used to compare with the adjoint method fairly. Using operator models attached in the supplemental material, the results in this paper can be reproduced. 
\subsection{Source control of the Poisson equation}\label{appendix_poisson}

During Phase 1, we generate the control input $m$ using the Fourier basis with frequency at most 3.
From this $m$, we generate the numerical label $u(m)$ using the finite difference method. It consists of 1000 training data and 200 test data with a resolution of $64\times64$. During this experiment, we use both the supervised loss \eqref{loss_supervised} and the residual loss \eqref{loss_residual} to compare them. For the residual loss, the laplacian $\Delta$ in the Poisson equation \eqref{poisson_equation} is approximated by a convolutional layer with the kernel
\begin{align*}
  \frac{1}{h^2}
  \begin{bmatrix}
    0 & 1 & 0\\
    1 & -4 & 1\\
    0 & 1 & 0
  \end{bmatrix},
\end{align*}
where $h$ is the mesh size. In addition, to force $u_\theta(m)$ to satisfy the Dirichlet boundary condition, the discrete metric between each pixel and the boundary is multiplied pixelwise \cite{berg2018unified}.

We use a encoder $G_\theta^{enc}$ which consists of 4 convolutional layers and a linear layer. Two decoders $G_\theta^{sol}$ and $G_\theta^{rec}$ consist of 4 convolution layers and 2 upsampling layers, respectively. To optimize paramter $\theta$, we use Adam optimzier with learning rate $10^{-3}$ and weight decay $10^{-6}$. A scheduler is used with step size 300 and discount factor 0.5. For both loss functions \eqref{loss_supervised} and \eqref{loss_residual}, we use $\lambda_1 = 1.5$. As shown in Figure \ref{poisson_operator}, the solution operators are well approximated for both cases, data-driven and data-free.

During Phase 2, the learned parameter $\theta^*$ are fixed and the control input $m$ is optimized. We compare the three initial guesses $m_{init} = xy(1-x)(1-y), x+y$, and $xy$. The L-BFGS optimizer is used with learning rate 1 and the regularizer weight $\lambda_2$ is 0.005.

One interesting result is that our method is robust to the initial guess for control input $m$ compared to the adjoint method. Table \ref{poisson_table} shows the corresponding results. It gives the control results for three different initial guesses. Our method gives similar results regardless of the initial guess. On the other hand, the relative error between the obtained optimal control $m$ and $m^*$ varies greatly when using the adjoint method. In the case $m_{init} = xy(1-x)(1-y)$, the adjoint method shows the most accurate prediction, but in the other cases, it predicts very far from the analytic optimal $m^*$ and is slower than our method.

\begin{figure*}[t]
\centering
\includegraphics[width=0.9\textwidth]{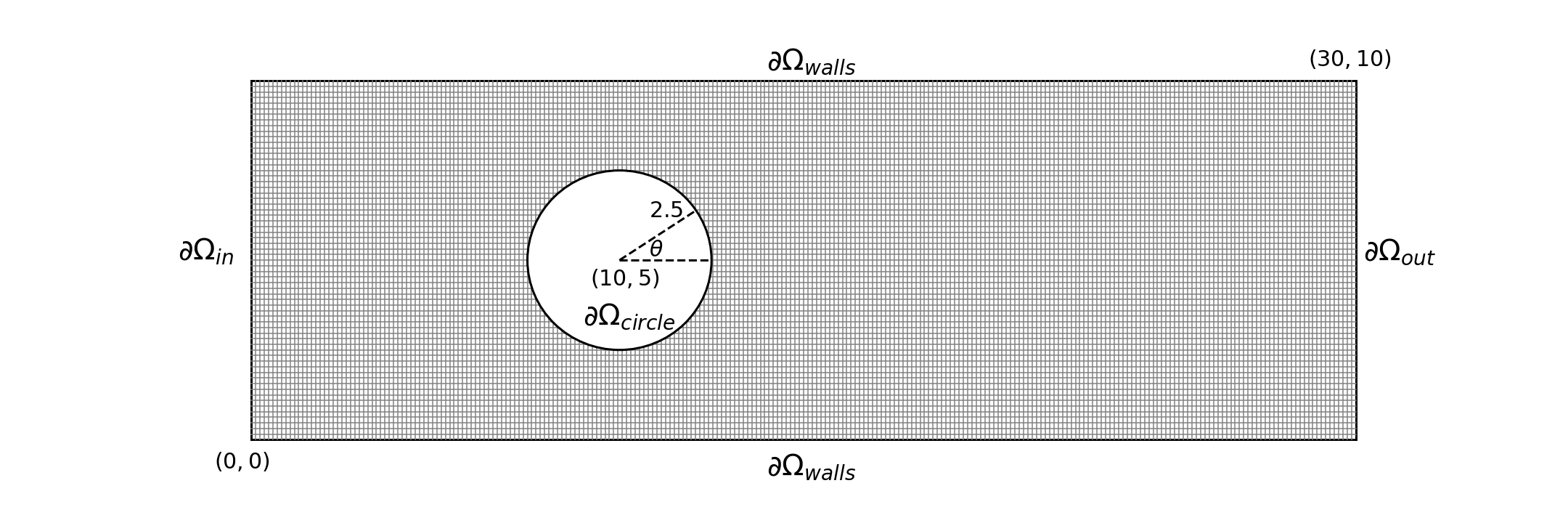} 
  \caption{The domain configuration of the boundary control of the Stokes equation \eqref{stokes_governing}.}
  \label{stokes_domain}
\end{figure*}

\begin{figure*}[t]
\centering
\includegraphics[width=0.8\textwidth]{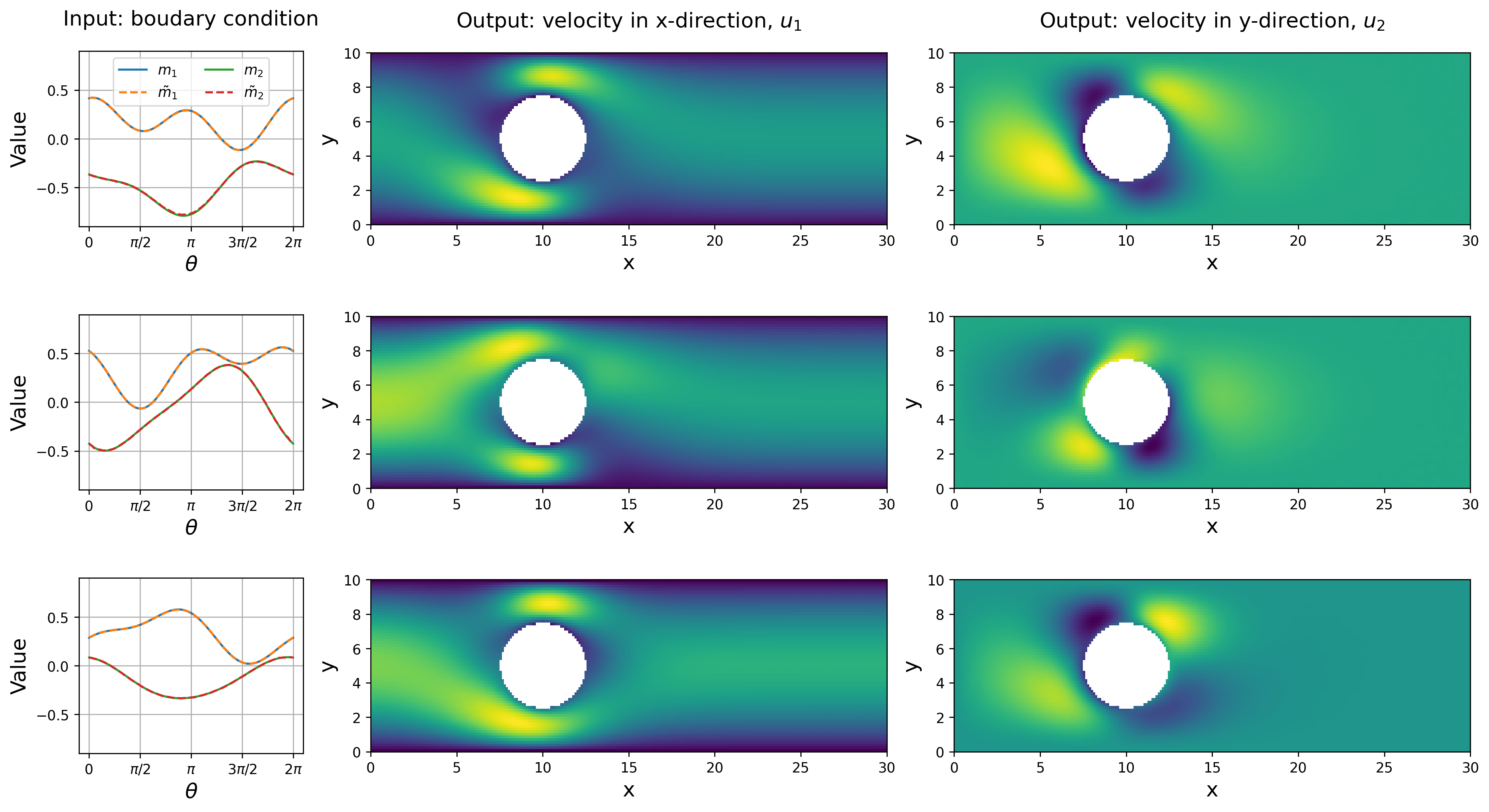}
  \caption{Results of operator learning for the Stokes equation. The first column is the control input (boundary condition on the circle), and the second and third columns are the corresponding output state.}
  \label{stokes_operator}
\end{figure*}

\begin{figure}[t]
\centering
\includegraphics[width=0.7\columnwidth]{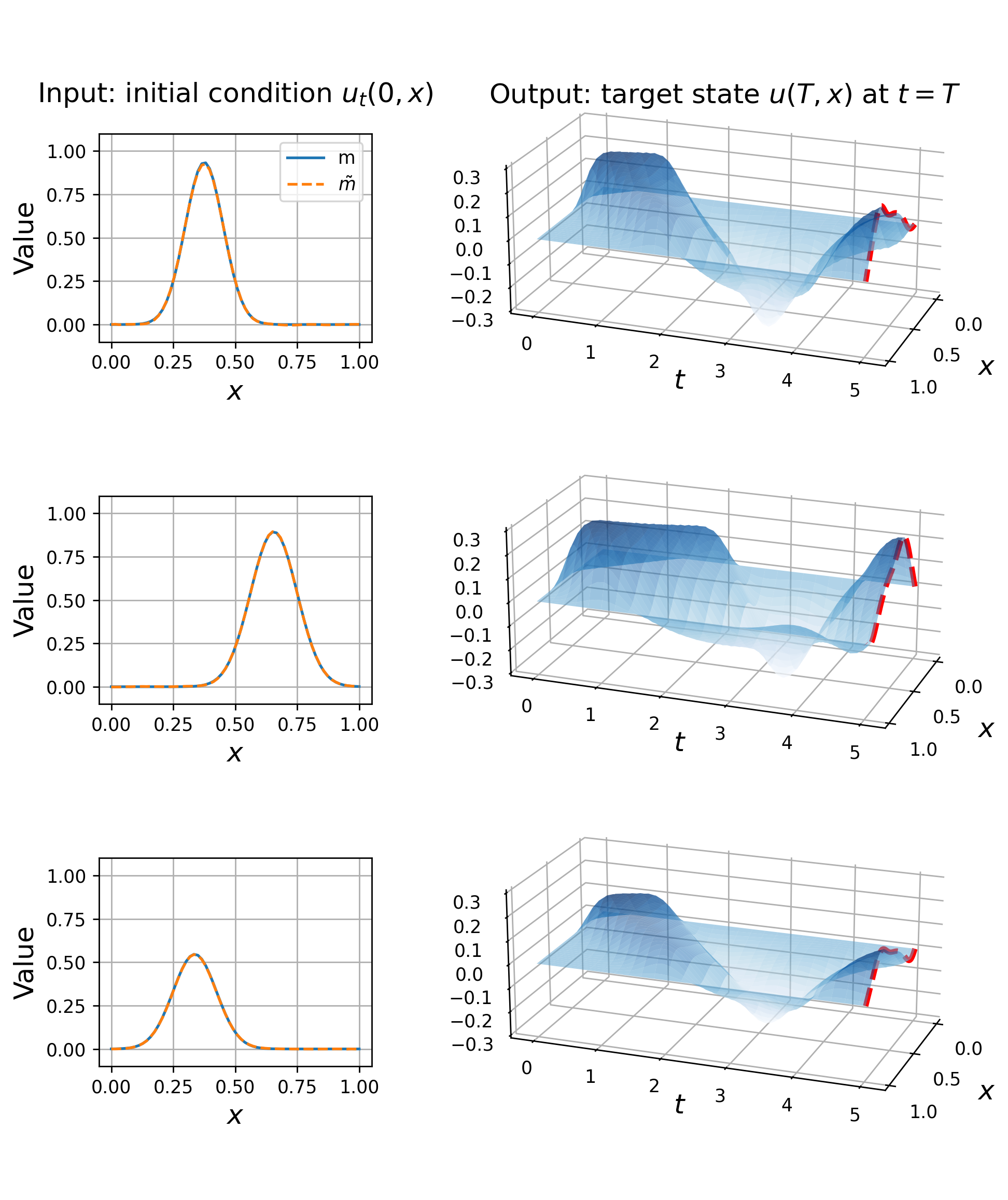} 
  \caption{Results of operator learning for the wave equation. The first column is the control input (initial condition $u_t(0,x)$), and the second column is the corresponding output state. The red dashed lines indicate the prediction of $u(T, x)$.}
  \label{wave_operator}
\end{figure}

\subsection{Boundary control of the stationary Stokes equation}\label{appendix_stokes}

To generate the boundary value on the circle $\Omega_{circle}$, the control input $m: \mathbb{R}^2 \rightarrow \mathbb{R}^2$ are transformed into a function of $\theta$, where $\theta$ is an angle from the $x$-axis centered at $\Omega_{circle}$ (See Figure \ref{stokes_domain} for the domain configuration). We sampled 1200 data for $m(\theta)$ of the form $m(\theta) = \sum_{n=1}^{d} c_n \sin^n \theta + d_n \cos^n \theta$ with $c_n, d_n \sim U([-0.3, 0.3])$ and $d=2$. The label is obtained using the finite element method.

The surrogate model for the solution operator consists of one encoder $G_\theta^{enc}$ and two decoders $G_\theta^{sol}$ and $G_\theta^{rec}$. The networks $G_\theta^{enc}$ and $G_\theta^{rec}$ consist of 1d-convolution layers, of which the input and output layer have two channels. Each channel contributes to the $x$ and $y$ directions of the control input vectors. The network $G_\theta^{sol}$ has four 2d-convolution layers, of which the output layer has two channels for representing the velocity in $x$ and $y$ directions. 

During Phase 1, the Adam optimizer with learning rate $10^{-3}$ is used with a scheduler. The learning rate is halved for every 100 epochs. We used three different scale meshes, namely $h \in [64, 128, 256]$ are the resolution of each mesh. Figure \ref{stokes_operator} shows the operator learning results for test data (when $h=256$). The first column describes the boundary value on the circle, where the solid line corresponds to the control input and the dashed line corresponds to the reconstruction by our network. The second and third columns are the visualization of the output velocity. The relative errors for test data are 0.012, 0.006, and 0.004, respectively for different values of $h$.

During Phase 2, the optimal boundary condition is recovered using the surrogate model trained during Phase 1. The L-BFGS optimizer with learning rate $10^{-2}$ is used to solve the inverse problem. The strong Wolfe line search is also adopted for the optimization. The hyperparameter $\lambda_2$ is chosen to be 10. The results for Phase 2 are explained in Section \ref{experiment}. Figure \ref{stokes_control} summarize the quantitative comparison of our method to the adjoint method. This highligted that our method achieve much better computaional efficiency compared to the adjoint method.

\subsection{Inverse design of nonlinear wave equation}\label{appendix_wave}

To train the surrogate model for the solution operator of the nonlinear wave equation, we collect the control-state data pairs consisting of the different initial conditions. The initial displacement $u(x,0)$ is fixed to zero, but the initial velocity $u_t(x,0)$ varies for each data in the form of a Gaussian pulse. We first sample 3000 data of $u_t(x,0)$, which are shaped as a Gaussian distribution with different parameters. The labels are then obtained through the finite element methods with the central difference approximations of time.

The surrogate model consists of one encoder $G_\theta^{enc}$ and two decoders $G_\theta^{sol}$ and $G_\theta^{rec}$. Each network has four 1d-convolution layers. The Adam optimizer with the learning rate $10^{-3}$ and the weight decay $10^{-4}$ is used during Phase 1. The learning rate scheduler is also used. Figure \ref{wave_operator} shows the samples of the solution operator predictions for the wave equation. The first column describes the control input, which corresponds to initial condition for $u_t(0, x)$, and the second column is the trajectory of state for each initial condition. The intermediate states are recovered using the numerical solver. Our model evaluate the state at target time $t=T$ much faster than the numerical solver, since it does not require time evolution with small time step. As shown in the figure, the reconstruction and the solution approximation work well, leading to the successful optimal control search. The relative error for the test data is 0.0068.

During Phase 2, we fix the surrogate model from Phase 1, and search optimal initial velocity $u_t(x, 0)$. The L-BFGS optimizer with strong Wolfe line search algorithm is adopted. The hyperparameter for the loss is chosen as $\lambda_2 = 0.001$. For comparison, the adjoint-based numerical approach is implemented. The most time-consuming part of the adjoint method is solving a nonlinear variational problem every time step. Newton iteration convergence is required for every iteration until it reaches the target time $T=5$. This significantly increases the inference time when using the numerical adjoint method. See Table \ref{wave_table} for summary of the results.

\begin{figure*}[t]
\centering
\includegraphics[width=0.8\textwidth]{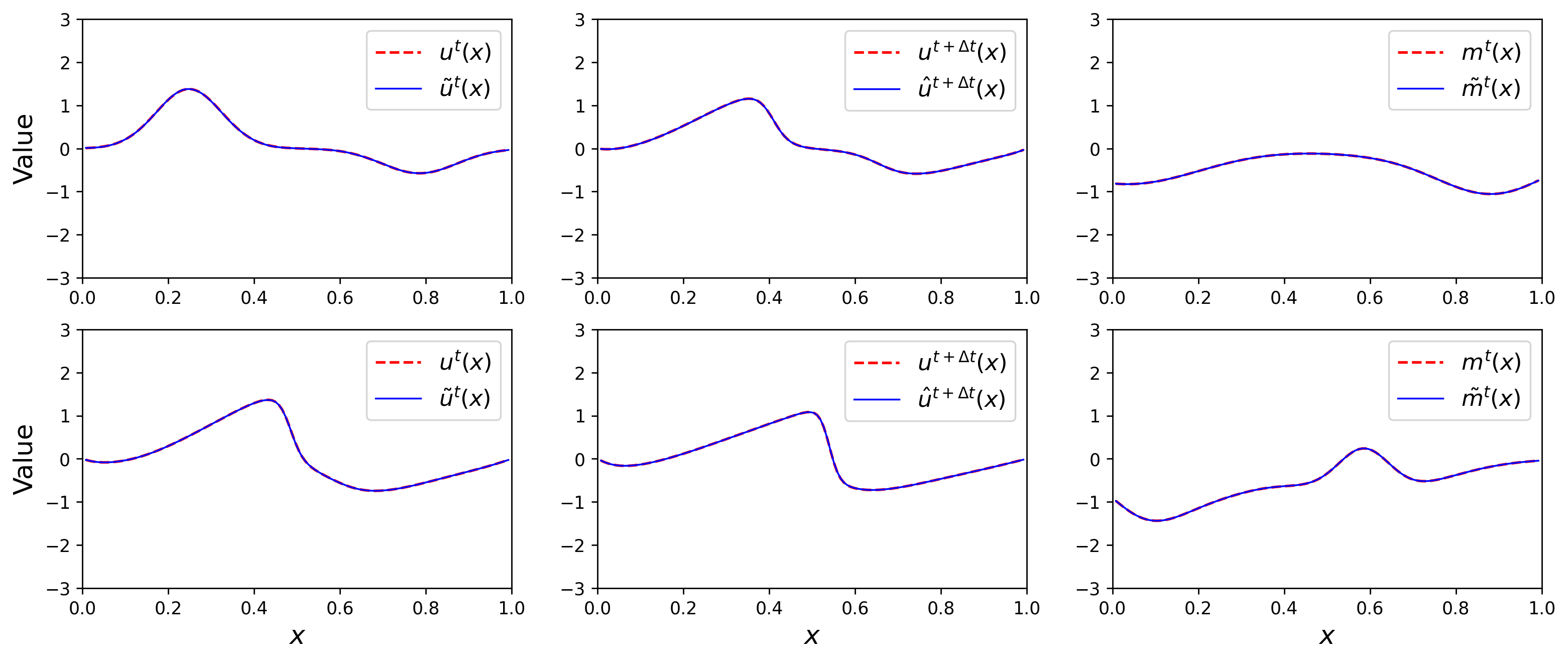} 
  \caption{Two samples for the results of operator learning during Phase1 to Burgers' equation}
  \label{burgers_operator}
\end{figure*}

\begin{figure*}[t]
\centering
\includegraphics[width=0.8\textwidth]{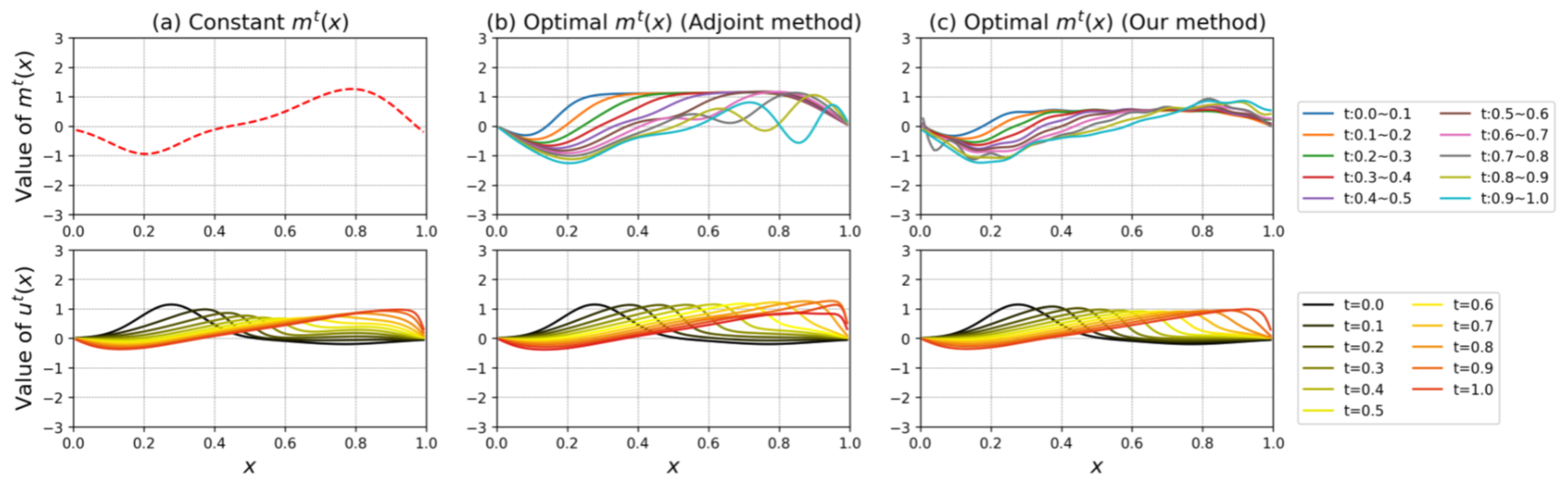} 
  \caption{The external force and trajectories with the external force using Burgers' equation for different target $u_d(x)$.}
  \label{burgers_control_ex2}
\end{figure*}

\subsection{Force control of Burgers’ equation}\label{appendix_burgers}

To train the surrogate model for approximating the solution operator during Phase 1, we use a dataset $\{(u^t,m^t,u^{t+\Delta t})\}$, which contains state tuples with the corresponding external force (control input). The state $u^t$ moves to the state $u^{t+\Delta t}$ by applying the external force $m^t$. The dataset is generated using the finite difference method with 130 resolutions containing two boundary points. Since the Dirichlet boundary condition is enforced at every time step, the 128 resolution on the interior points is used for both the state $u^t$ and the force $m^t$. For the initial state $u^0$, two Gaussian waves are added, which have a center at $+x$ with negative amplitude and at $-x$ with positive amplitude to generate a shock. In addition, the external force $m^t$ is randomly generated by adding 1$\sim$7 Gaussian distributions randomly. Using 60 random initial conditions and 500 external forces, the 30000 trajectories are generated to construct a dataset. The force $m^t$ is fixed for each initial condition when the trajectory is generated. The first column plots in Figure \ref{burgers_control_ex1} show an example of the fixed force $m^t$ (upper plot) and the corresponding trajectory of the state $u^t$.

During Phase 1, both the two encoders ($H_{\theta_1}^{enc}$, $G_{\theta_2}^{enc}$) and the two decoders ($H_{\phi_1}^{rec}$, $G_{\phi_2}^{rec}$) have four 1d-convolution layers each. The transition network $T_{\psi}^{tran}$ has one 1d-convolution layer. The Adam optimizer with the learning rate $10^{-3}$ and the weight decay $10^{-4}$ with the scheduler is used for the solution operator learning in Phase 1. The weights for each term in the loss function \eqref{loss_burgers_phase1} are chosen as 1, 1.5, 0.5, 1.

Figure \ref{burgers_operator} shows the results of the trained solution operator for Burgers' equation. We plot the two samples of the dataset in each row. The first and third columns in Figure \ref{burgers_operator} show the state $u^t$ and the external force $m^t$ (blue line), the pairing of which is in our input dataset and along with their reconstruction (red dotted line), which are the outputs of both autoencoders for the state and external force. As the figure shows, the reconstructions of the state $u^t$ and of the force $m_t$ are well approximated. The second column plots in Figure \ref{burgers_operator} show the approximated state of the next time step $u^{t+\Delta t}$ and its reconstruction. It shows that the transition of the state $u^t$ to the next step is well approximated using the extended model indicated in Figure \ref{architecture_time}.

During Phase 2, the L-BFGS optimizer with strong Wolfe line search algorithm, the learning rate $0.5$ and the tolerance change $4\times10^{-6}$ is used to search the optimal external force. The weights for each term in the loss function \eqref{loss_burgers_phase2} are chosen as 10 and 0.3. We compare the results of optimal force inferred by our trained model to the adjoint method. The tolerance change and other settings of the L-BFGS optimizer were the same for both the our model and the numerical method. Figure \ref{burgers_control_ex1} and Figure \ref{burgers_control_ex2} show the results of optimal external force and trajectory for different targets $u_d(x)$ compared to the results of the adjoint method. Table \ref{burgers_table} summarizes the quantitative comparison of our method to the adjoint method. 

\end{document}